\documentclass{agn_article}

\usepackage{graphicx}
\usepackage{bm}
\usepackage{bbold}
\usepackage{subfigure}
\usepackage{float}
\usepackage{amsfonts,color}
\usepackage{epsfig}
\usepackage{lscape}
\usepackage{footnote}
\usepackage[utf8]{inputenc}	
\usepackage{mathtools}
\usepackage{wasysym}
\usepackage[english]{babel}
\usepackage{multirow}
\usepackage{array}
\usepackage{colonequals}
\usepackage{authblk}
\usepackage{booktabs}

\newcommand{\tablebreak}{\\[-1em]\\[.1em]}

\newcolumntype{C}[1]{>{\centering\let\newline\\\arraybackslash\hspace{0pt}}m{#1}}

\DeclareMathOperator{\dev}{dev}
\def\tr{\textrm{tr}}
\def\dd{\displaystyle}

\begin{document}

\title{The exponentiated Hencky strain energy in modelling tire derived 
material for moderately large deformations\footnote{To
appear in ASME Journal of Engineering Materials and
Technology}}

\author{Giuseppe Montella\thanks{Department of Structure for Engineering and Architecture, University of Naples `Federico II', Naples, Italy;\;\\ Department of Civil and Environmental Engineering, University of California Berkeley, Berkeley, CA,USA;\;\\ email: giuseppe.montella@unina.it},\quad
Sanjay Govindjee\thanks{Department of Civil and Environmental Engineering, University of California Berkeley, Berkeley, CA, USA;\; ASME member, email: s\_g@berkeley.edu}\quad and\quad
Patrizio Neff\thanks{Faculty of Mathematics, University of Duisburg-Essen, Essen, Germany;\; email:
patrizio.neff@uni-due.de}}
\date{\today}

\maketitle

\begin{abstract}

This work presents a hyper-viscoelastic model, based on the Hencky-logarithmic strain tensor to model the response of a Tire Derived Material (TDM) undergoing moderately large deformations. TDM is a composite made by cold forging a mix of rubber fibers and grains, obtained by grinding scrap tires, and polyurethane binder. The mechanical properties are highly influenced by the presence of voids associated with the granular composition and low tensile strength due to the weak connection at the grain-matrix interface. For these reasons, TDM use is restricted to applications concerning a limited range of deformations. Experimental tests show that a central feature of the response is connected to highly nonlinear behavior of the material under volumetric deformation which conventional hyperelastic models fail in predicting. The strain energy function presented here is a variant of the exponentiated Hencky strain energy, which for moderate strains is as good as the quadratic Hencky model and in the large strain region improves several important features from a mathematical point of view. The proposed form of the exponentiated Hencky energy possesses a set of parameters uniquely determined in the infinitesimal strain regime and an orthogonal set of parameters to determine the nonlinear response.
The hyperelastic model is additionally incorporated in a finite deformation viscoelasticity framework that accounts for the two main dissipation mechanisms in TDMs, one at the microscale level and one at the macroscale level. The new model is capable of predicting different deformation modes in a certain range of frequency and amplitude with a unique set of parameters with most of them having a clear physical meaning. This translates into an important advantage with respect to overcoming the difficulties related to finding a unique set of optimal material parameters as are usually encountered fitting polynomial forms of strain energies. Moreover, by comparing the predictions from the proposed constitutive model with experimental data we conclude that the new constitutive model gives accurate prediction.
\end{abstract}

\newpage

\section{Introduction}

\label{sec:intro}
\begin{figure}
        \centering
        \subfigure[\,TDM pad]{\includegraphics[height=3cm]{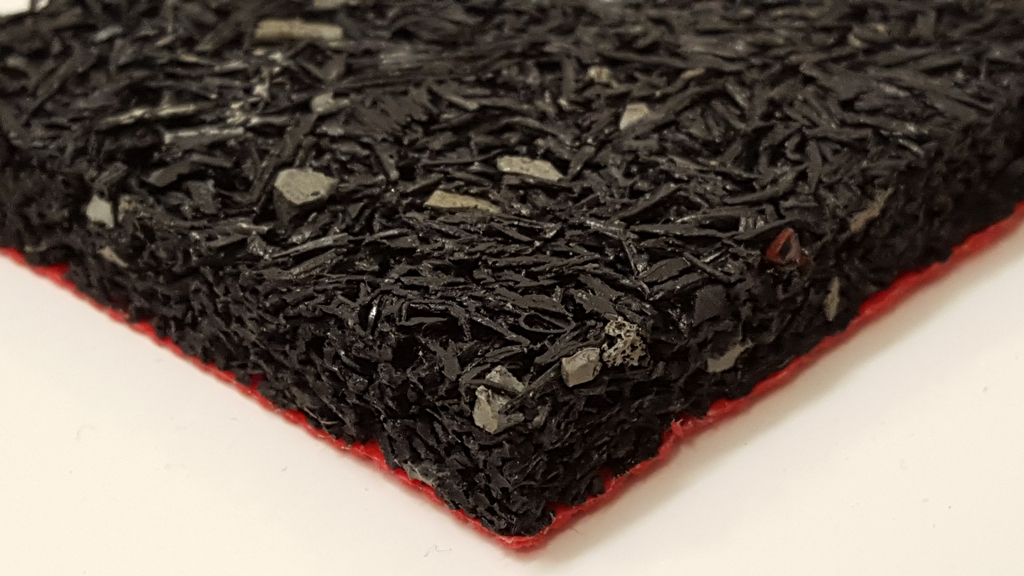}}\hspace{2cm}
        \subfigure[\,Close up view]{\includegraphics[height=3cm]{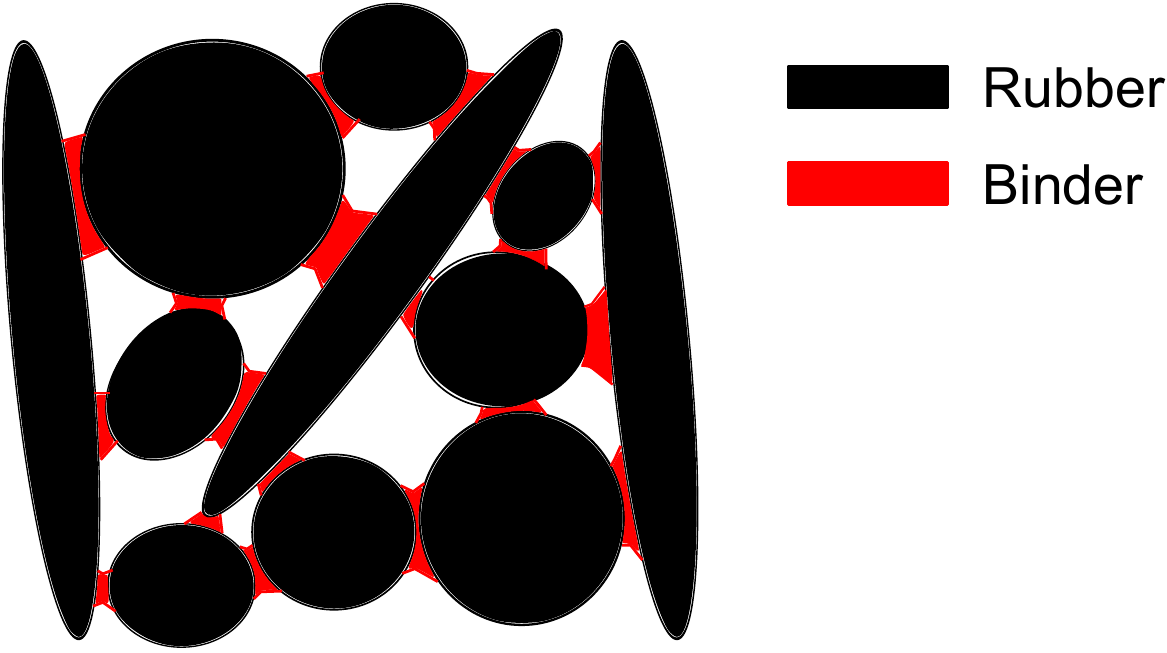}}
        \caption[Tire Derived Material]{Tire Derived Material.}\label{fig:fig_TDM} 
        \end{figure}

In spite of a rapid growth in technologies and development, scrap tire disposal is still an important and unresolved environmental engineering issue today. One promising component, from the spectrum of proposed solutions, is the recycling of tires into engineering materials. In this paper, the mechanical characterization of a new Tire Derived Material (TDM) for structural applications is proposed. This TDM is obtained by grinding scrap tires and rubber factory leftovers into grains and fibers, these together with a polyurethane binder are first leveled by a roller and then pressed together to form TDM pads (Figure \ref{fig:fig_TDM}). The material, made of Styrene-Butadiene Rubber (SBR), the most popular rubber in tire production, has a low cost and easy to implement production cycle -- alternately Ethylene-Propylene Diene Monomer (EPDM), a rubber used for the production of a wide variety of seals can be used. TDMs can be made in various densities with different mechanical properties and have been used mainly in railway applications for vibration reduction \cite{Montella2012}.
The usual composition of TDMs results in high compressibility and allows for use only in a moderately large range of deformation.
  
Unfortunately the common hyperelastic material models, e.g.~Arruda-Boyce \cite{Arruda1993}, Mooney-Rivlin \cite{mooney1940}, and Ogden \cite{ogden1972} models, fail in describing their behavior in different deformation modes with a unique set of parameters \cite{Montella2014}; see also Appendix \ref{sec:appA}. It is also noted that fitting experimental data of elastomeric solids to polynomial-like strain energy functions is not an easy task and can lead to oscillating functions with parameters that may not have physical meaning \cite{Ogden2004}. Here, a logarithmic measure is used to describe the mechanical behavior of TDMs.

Logarithmic strain, typically referred as ``true strain", was first applied to elasticity theory by the geologist G.F.~Becker \cite{becker1893,neff2014becker}, who was an instructor of mining and metallurgy at Berkeley from 1875 to 1879. However, its introduction is often attributed to P.~Ludwik \cite{Ludvik}, who defined (one-dimensional) logarithmic strain via the integral $\int_{l_{0}}^l \! \frac{dl}{l}$ in order to measure the extension of a rod of length \textit{l}. Today, the logarithmic strain tensor is also named after H.~Hencky \cite{Hencky_1928,Hencky1929}, who used it in his systematic deduction of an idealized elastic law \cite{henckyTranslation}. The Hencky strain measure has many interesting properties, one of the most useful is that it allows for the full realization of an uncoupled additive split of volumetric and deviatoric deformations at finite strain. The elastic law proposed by Hencky, which is in good agreement with experiments for a wide class of materials for moderately large deformations, as Anand demonstrated \cite{Anand_79}, is induced by the so-called  quadratic Hencky strain energy:
    \begin{equation}\label{hencky}
    \ W_{\mathrm{H}}(F)\colonequals\widehat{W}_{\mathrm{H}}(U)\colonequals\mu\,\|\dev_3 \log {U}\|^{2} +\frac{\kappa}{2}\,[\mathrm{tr}(\log {U})]^{2}.  
    \end{equation}
    
As observed by Hsu, Davies and Royles \cite{hsu.ea:67}, the choice of strain measure can facilitate the transference of the nonlinearity in the stress-strain response from the strain-to-stress mapping to the strain measure itself; see also Sharda \cite{sharda:74}.  This in part explains the success of the quadratic Hencky strain energy approach.  Notwithstanding, the Hencky strain energy also possesses interesting intrinsic mathematical properties.
In a series of articles \cite{NeffGhibaLankeit,NeffGhibaPoly,NeffGhibaPlasticity,NeffGhibaAdd}, a family of isotropic volumetric-isochoric decoupled strain energies
    \begin{align}\label{the}
    W_{\mathrm{eH}}(F)\colonequals\widehat{W}_{\mathrm{eH}}(U)\colonequals\dd\left\{\begin{array}{lll}
    \dd\frac{\mu}{k}\,e^{k\,\|\dev_n\log {U}\|^2}+\dd\frac{\kappa}{{2\, {\widehat{k}}}}\,e^{\widehat{k}\,[\tr(\log U)]^2}&\text{if}& \det\, F>0\vspace{2mm}\\
    +\infty &\text{if} &\det F\leq 0
    \end{array}\right.\quad
    \end{align}
based on the Hencky-logarithmic strain tensor $\log U$ were studied. Here $\mu>0$ is the infinitesimal shear modulus,
$\kappa=\frac{2\mu+3\lambda}{3}>0$ is the infinitesimal bulk modulus with $\lambda$ the first Lam\'{e} constant, $k,\widehat{k}$ are dimensionless
parameters, $F=\nabla \varphi$ is the gradient of deformation,  $U=\sqrt{F^T F}$ is the right stretch tensor and $\dev_n\log {U} =\log {U}-\frac{1}{n}
 \tr(\log {U})\cdot\mathbb{1}$
is the deviatoric part  of the logarithmic\footnote{Here and throughout, $\log$ denotes the natural logarithm.} strain tensor $\log {U}$. This family of exponentiated Hencky strain energies improves upon the well-known properties of the original Hencky strain energy.
In particular, it was recently found that the Hencky energy (not the logarithmic strain itself) exhibits  a fundamental property: by purely differential geometric reasoning, it was shown \cite{NeffEidelOsterbrinkMartin_Riemannianapproach,Neff_Osterbrink_Martin_hencky13,Neff_Nagatsukasa_logpolar13} (see also  \cite{Neff_log_inequality13,LankeitNeffNakatsukasa,borisov2015}) that
    \begin{align}\label{geoprop}
    \mathrm{dist}^2_{\mathrm{geod}}\left((\det F)^{1/n}\cdot \mathbb{1},  \mathrm{SO}(n)\right)&= \mathrm{dist}^2_{ \mathrm{geod,\mathbb{R}_+\cdot \mathbb{1}}}\left((\det F)^{1/n}\cdot \mathbb{1}, \mathbb{1}\right)= \frac1n\, [\mathrm{tr}(\log U)]^2 = \frac1n\, (\log (\det U))^2,\notag\\
    \mathrm{dist}^2_{\mathrm{geod}}\left( \frac{F}{(\det F)^{1/n}}, \mathrm{SO}(n)\right)&=\mathrm{dist}^2_{ \mathrm{geod}, \mathrm{SL}(n)}\left( \frac{F}{(\det F)^{1/n}}, \mathrm{SO}(n)\right)=\|\dev_n \log U\|^2,
    \end{align}
where $\mathrm{dist}_{\mathrm{geod}}$ is the canonical left invariant geodesic distance on the Lie group $\mathrm{GL}^+(n)$ and $\mathrm{dist}_{\mathrm{geod},\mathrm{SL}(n)}$, $\mathrm{dist}_{\mathrm{geod},\mathbb{R}_+\cdot \mathbb{1}}$ denote the corresponding geodesic distances on the Lie groups $\mathrm{SL}(n)$ and $\mathbb{R}_+\cdot \mathbb{1}$, respectively (see \cite{Neff_Osterbrink_Martin_hencky13,Neff_Nagatsukasa_logpolar13}). Thus $W_{\mathrm{H}}$ and $W_\mathrm{eH}$ have the attractive feature that the energies are based directly on a geometrically intrinsic distance of the deformation gradient to the group of rigid rotations.

For small elastic strains, $W_{\mathrm{eH}}$ approximates the classical quadratic Hencky strain energy $W_{\mathrm{H}}$, which is not everywhere rank-one convex; moreover in \cite{NeffGhibaLankeit}, it is also pointed out that the quadratic Hencky energy has some other serious shortcomings.
These points being more or less well-known, it is clear that there cannot exist a general mathematical well-posedness result for the quadratic Hencky model $W_{\mathrm{H}}$, although an existence proof for small loads based on the implicit function theorem is, of course, possible. The use of \eqref{the} allows for the retention of the fundamental geometric property \eqref{geoprop} of the original Hencky strain energy, but at the same time alleviates some of its mathematical drawbacks: up to moderate strains, {for principal stretches $\lambda_i\in(0.7,1.4)$}, the exponentiated Hencky formulation \eqref{the} is \emph{de facto} as good as the quadratic  Hencky model $W_{\mathrm{H}}$, and in the large strain region it improves several important features from a mathematical point of view. The main feature is that the exponentiated Hencky energy \eqref{the} satisfies the Legendre-Hadamard condition (rank-one convexity) in planar elasto-statics \cite{NeffGhibaLankeit,isochoricRankOne}, i.e.\ for $n=2$. In this case, the energy is even polyconvex, which, together with a coercivity estimate, allows for the application of classical theorems for the existence of energy minimizers \cite{NeffGhibaPoly,neffGhibaSilhavy2015improvement}.

Despite these advantages, some aspects of the three-dimensional description remain open, since the formulation is not globally rank-one convex. However, in the three-dimensional case, a loss of ellipticity only occurs for extreme distortional strains \cite{ghibaEllipticityDomain}. This suggests that the exponentiated Hencky energy \eqref{the} retains its full suitability for materials that undergo additional (typically irreversible) phenomena based on distortional criteria of the Huber-Hencky-von-Mises type, as the involved elasticity tensors can thereby be prevented from reaching the non-elliptic domain. This is in sharp contrast to the loss of ellipticity of the quadratic Hencky energy $W_{\mathrm{H}}$, which is not related to the distortional energy alone. 

Beside the above mathematical advantages, the exponentiated Hencky energy satisfies a number of additional desirable constitutive properties \cite{NeffGhibaLankeit}: for example, planar pure Cauchy shear stress always induces biaxial pure shear strain; the limit case $\kappa \rightarrow +\infty$ or, equivalently, $\nu=\frac{1}{2}$ for the linear Poisson's ratio $\nu$, corresponds to exact finite incompressibility; and there exists a certain three parameter subset ($k=\frac{2}{3}\, \widehat{k}\,$) such that uniaxial tension leads to no lateral contraction if and only if $\nu=0$ (i.e.\ $\kappa = \frac{2}{3}\mu$), as in linear elasticity (see [24] for further discussion). Like the quadratic Hencky energy, the exponentiated Hencky energy also satisfies a weakened version of Truesdell's empirical inequalities \cite{borisov2015}.

In this paper, a variation to the volumetric part of \eqref{the} is proposed to capture the high nonlinearity of TDMs when subjected to volumetric deformation:
    \begin{align}\label{them}
    W_{ \mathrm{eHm}}(F)\colonequals\widehat{W}_{ \mathrm{eHm}}(U)\colonequals\dd\left\{\begin{array}{lll}
    \dfrac{\mu}{k}\, e^{k\,\|\dev_n\log {U}\|^2}+\dfrac{\kappa}{2\,\widehat{k}}\,e^{\widehat{k}\,[\mathrm{tr}(\log  {U})]^{2}}+\dfrac{\kappa_{1}}{m\,\tilde{k}}\,e^{\tilde{k}\, |\mathrm{tr}(\log {U})|^{m}}&\text{if} &\det\, F>0,\vspace{2mm}\\
    +\infty &\text{if} &\det F\leq 0,
    \end{array}\right.
    \end{align}
where $\kappa_{1}$ is the value of the bulk modulus for large deformations and $m$ and $\tilde{k}$ are dimensionless parameters. The main advantage of using the modified exponentiated-Hencky energy comes from the fact that the shear and bulk modulus are already uniquely determined in the infinitesimal strain regime, while $\kappa_{1}$ determines the nonlinear response, without interfering with $\mu$ and $\kappa$. 
For the modified exponentiated-Hencky strain energy proposed in this work, the Kirchhoff stress tensor is given by:
    \begin{align}\label{tau}
    \tau= D_{\log {U}}{W}_{\mathrm{eHm}}(\log {U})= \,&2\,\mu \,e^{k\,\|\mathrm{dev}_{3}\log {U}\|^{2}}\cdot\dev_3\log {U} 
    \notag\vspace{2mm} \\&+\left[ \kappa\, e^{\widehat{k}\,[\mathrm{tr}(\log {U})]^{2}}\mathrm{tr}(\log {U}) 
    + \kappa_{1}\, e^{\tilde{k}\,|\mathrm{tr}(\log {U})|^{m}} \dfrac{|{\mathrm{tr}(\log {U})}|^{m}}{\mathrm{tr}(\log {U})}\right]  \cdot \mathbb{1}\, ,
    \end{align}
while the Cauchy stress tensor is:
    \begin{align}\label{sigma}
    \sigma= e^{\mathrm{-tr}(\log {U})}\tau= \,&2\,\mu \,e^{k\,\|\mathrm{dev}_{3}\log {U}\|^{2}\mathrm{-tr}(\log {U})}\cdot\dev_3\log {U}\notag\vspace{2mm} \\ &+
    \left[ \kappa\, e^{\widehat{k}\,[\mathrm{tr}(\log {U})]^{2}\mathrm{-tr}(\log {U})}\mathrm{tr}(\log {U}) 
    + \kappa_{1}\, e^{\tilde{k}\,|\mathrm{tr}(\log {U})|^{m}\mathrm{-tr}(\log {U})} \dfrac{|{\mathrm{tr}(\log {U})}|^{m}}{\mathrm{tr}(\log {U})}\right]  \cdot \mathbb{1}.
    \end{align}
    
\section{Rate-independent response}

In this section we discuss the rate-independent response of the material, which we will refer to as the equilibrium response. The physical properties of the TDMs are greatly influenced by the technologies used in manufacturing them. Tests have shown that the density and the mixture composition of the material are the parameters that most strongly affect its mechanical properties. In total we consider three types of TDMs with the same composition but different densities (Table \ref{tab:I}) and three modes of deformation: shear, uniaxial compression, and pseudo-hydrostatic compression.
The TDMs studied where made from tires that were shredded into chips, mostly 50 mm in size using a rotary shear shredder with two counter-rotating shafts;
tire chips were reduced to a size smaller than 10 mm in a granulator while most of the steel cords were liberated by a combination of shaking screens and wind shifters;
the styrene-butadiene rubber (SBR) granules were selected according to their dimensions to fit the desired design mix;
polyurethane binder was added to the rubber granules mix and uniformly distributed;
lastly pads of required size and shape were obtained by pressing the compound in a mold.
In the sections describing the tests, predictions of the model are also shown, all using a fixed set of fitted parameters which are discussed in Section \ref{sec:III}.
    \begin{table}[H]
    \caption{Tire Derived Material description.}\label{tab:I}
    \centering
    \begin{tabular}{C{2cm}C{2cm}C{3cm}}
    \toprule
    \multirow{2}{*}{Material} & Density & \multirow{2}{*}{Composition}\\
    &$\mathrm{(kg\cdot m^{-3})}$ &\tablebreak
    \hline\\[-.7em]
    TDM $500$ & $500$ &$90\%$ SBR fibers\tablebreak
    TDM $600$ & $600$ &+\tablebreak 
    TDM $800$ & $800$ &$10\%$ SBR grains\tablebreak
    \hline
    \end{tabular}
    \end{table}
    
\subsection{Simple shear}

Shear tests were carried out at \textit{Tun Abdul Razak Research Centre} (TARRC) in Hertford (UK). The samples were tested with the classical dual lap simple shear test configuration commonly used in the tire industry. Samples of 90 mm in width, 50 mm in length and 20 mm in thickness, were sheared to a shear strain amplitude of $100\% $ of the initial rubber thickness at the (slow) strain rate of $0.0067\, \mathrm{s^{-1}}$. The procedure is explained in detail in \cite{Montella2014}. In simple shear the direction of applied displacements does not coincide with the direction of the principal stretches; rather it involves a rotation of axes. The polar decomposition of $F=R \cdot U$ gives the right Biot stretch tensor $U=\sqrt{F^{T}F}$ of the deformation and the orthogonal polar factor $R$:
    \begin{equation}
    U=\frac{1}{\sqrt{\gamma^{2}+4}}
    \left[ {\begin{array}{ccc}
    2 & \gamma & 0 \\
    \gamma & \gamma^{2}+2 & 0\\
    0 & 0 & \sqrt{\gamma^{2}+4}
    \end{array} } \right] \,, \qquad
    R=\frac{1}{\sqrt{\gamma^{2}+4}}
    \left[ {\begin{array}{rcc}
    2 & \gamma & 0 \\
    -\gamma & 2 & 0\\
    0 & 0 & \sqrt{\gamma^{2}+4}
    \end{array} } \right].
    \end{equation}
$U$ can be orthogonally diagonalized to show:
    \begin{equation}
    \log {U}=\frac{1}{\sqrt{\gamma^{2}+4}}
    \left[ {\begin{array}{rrr}
     -\gamma  \log \lambda_{1} & 2\log \lambda_{1} & 0 \\
    2\log \lambda_{1} & \gamma \log \lambda_{1} & 0\\
     0 & 0 & 0
    \end{array} } \right],
    \end{equation}
where $\lambda_{1}=\dfrac{1}{2}\left (\sqrt{\gamma^{2}+4}+\gamma \right)$ is the first eigenvalue of $U$.
Simple shear does not involve a change in volume; for this reason $\det F=1$ and $\mathrm{tr}(\log {U})=0$. The non-zero Kirchhoff stress component $\tau_{12}$ from equation \eqref{tau} is given by:
    \begin{equation}\label{taushear}
    \tau_{12}=4\,\mu\, e^{2\,k\,\log^{2}\left[\dfrac{1}{2}\left (\sqrt{\gamma^{2}+4}+\gamma \right)\right]} \cdot  \frac{\log \left[\dfrac{1}{2}\left (\sqrt{\gamma^{2}+4}+\gamma \right)\right]}{\sqrt{\gamma^{2}+4}}.
    \end{equation}
     \begin{figure}[H]
        \subfigure[\,TDM 500]{\includegraphics{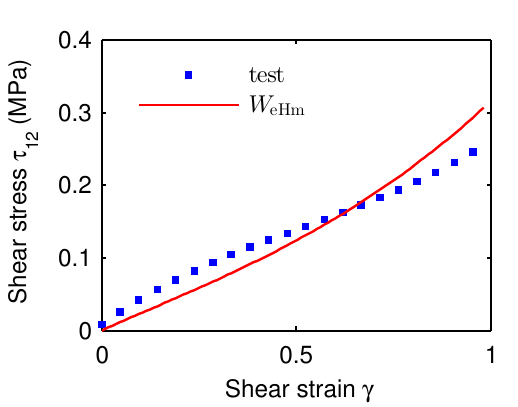}}
        \subfigure[\,TDM 600]{\includegraphics{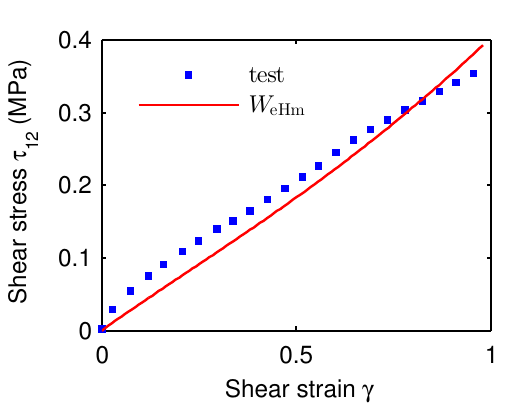}}
        \subfigure[\,TDM 800]{\includegraphics{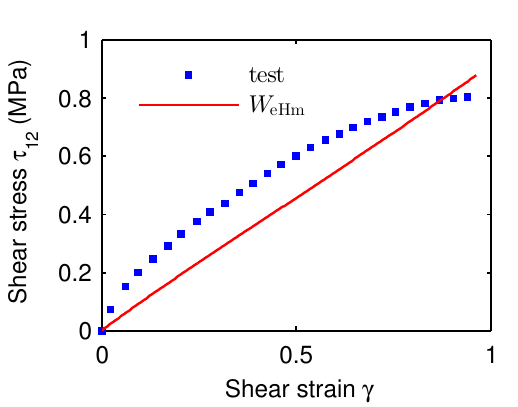}}
        \caption[Comparison between tests and exponentiated-hencky strain energy function]{Comparison between shear stress corresponding to exponentiated Hencky energy $W_{\mathrm{eHm}}$, equation \eqref{taushear}, and experimental tests for different densities.}\label{fig:fig_shear}
        \end{figure}
Figure \ref{fig:fig_shear} shows the ability of the model to capture the shear behavior out to a shear strain of 100\%. It is to be noted that the TDM 800 sample, Figure \ref{fig:fig_shear}(c), physically failed in the experiment due to cracking and crumbling.  Thus the poor correlation in Figure \ref{fig:fig_shear}(c) does not reflect poorly on the model. Unfortunately, an intact test is not available for TDM 800 in this configuration. Notwithstanding, given the reasonable agreement seen in Figures \ref{fig:fig_shear}(a) and \ref{fig:fig_shear}(b), we feel the model performs well in shear. This is in agreement with the findings in \cite{NeffGhibaLankeit} when $W_\mathrm{eH}$ was applied to the rubber data of \cite{treloar1944stress,jones1975properties}.

\subsection{Uniaxial compression}

\begin{figure}[H]
        \centering    
        \subfigure[\,Strain history]{\includegraphics{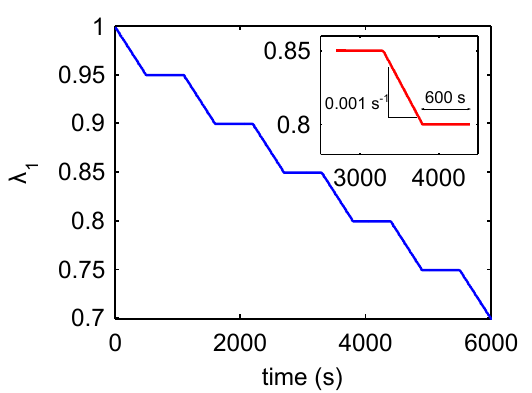}}
        \subfigure[\,TDM 500]{\includegraphics{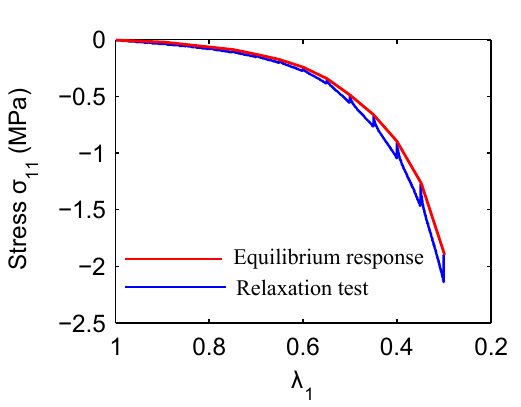}}\\
        \subfigure[\,TDM 600]{\includegraphics{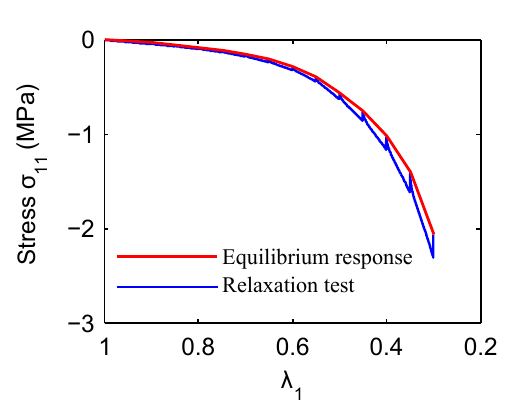}}
        \subfigure[\,TDM 800]{\includegraphics{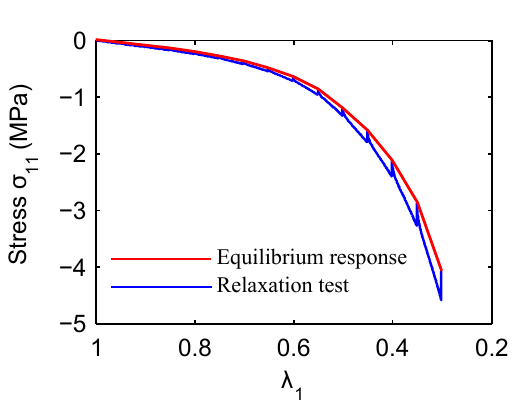}}
        \caption[Compression test]{Compression test procedure. First figure (a) shows strain history. Blue curve represent the true data.  Red curve represents the assumed equilibrium response from the data.} \label{fig:fig_load}
        \end{figure}

\begin{figure}[ht!]
        \centering     
        \subfigure{0\%}{\includegraphics[width=40mm,keepaspectratio]{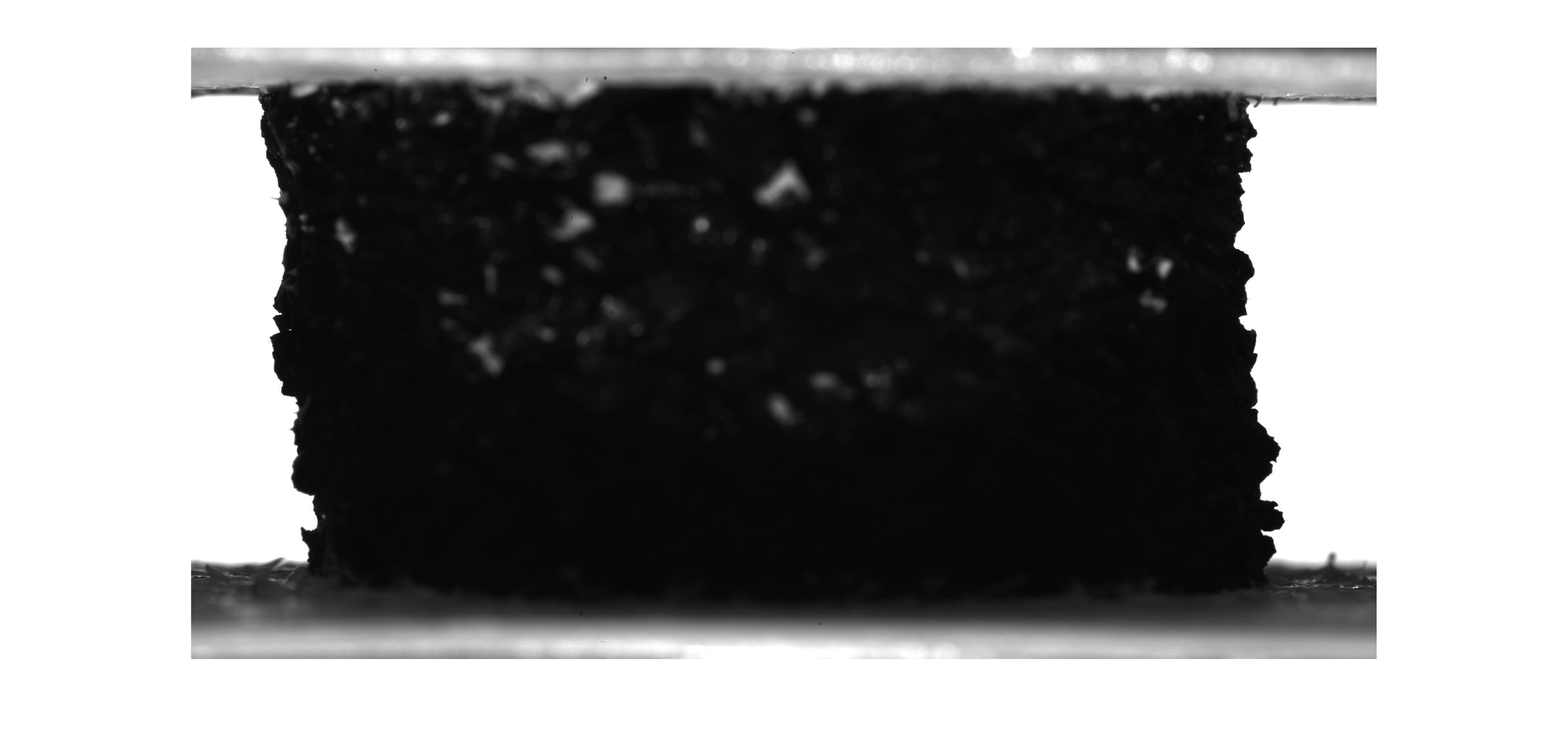}}
        \subfigure{\includegraphics[width=40mm,keepaspectratio]{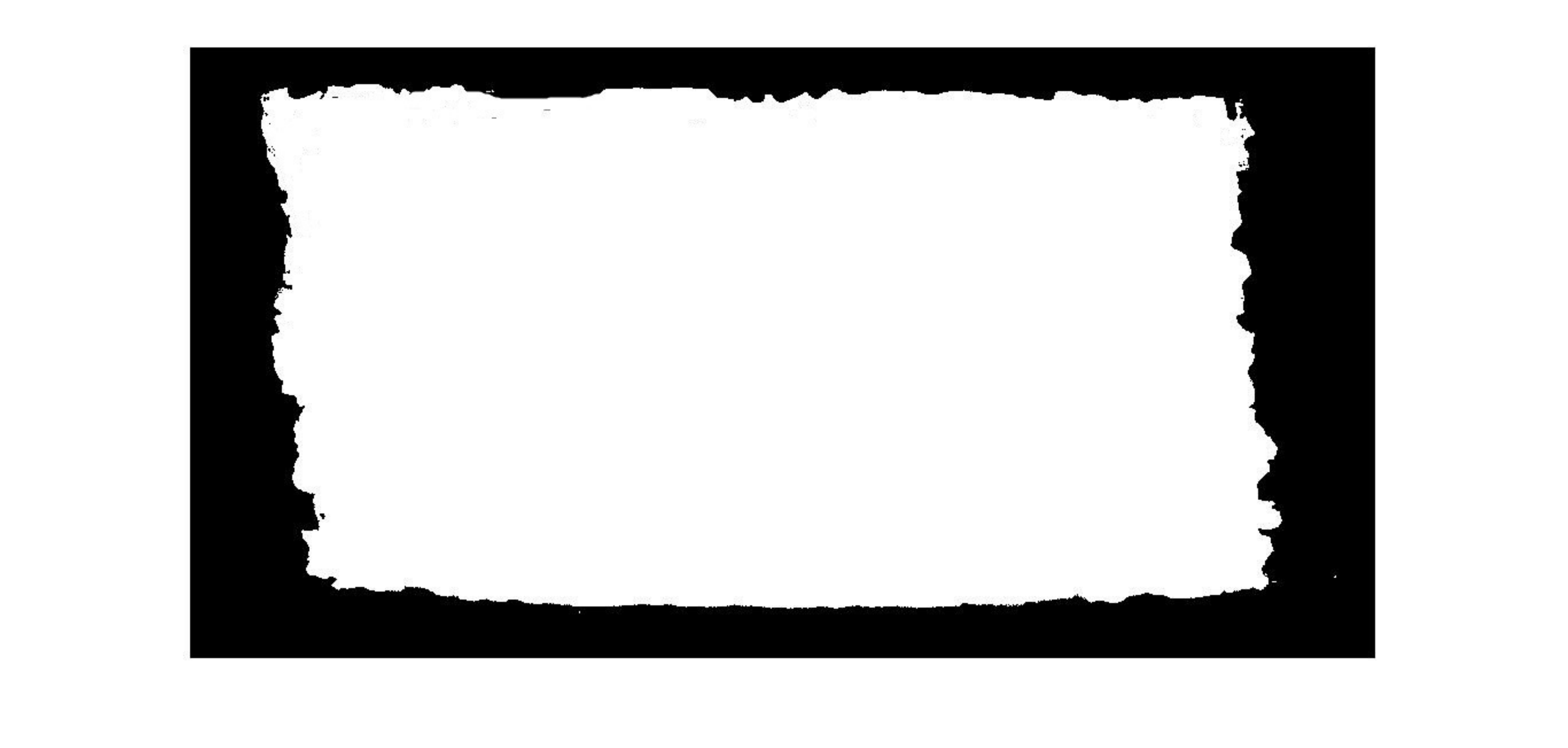}}\\
        \subfigure{35\%}{\includegraphics[width=40mm,keepaspectratio]{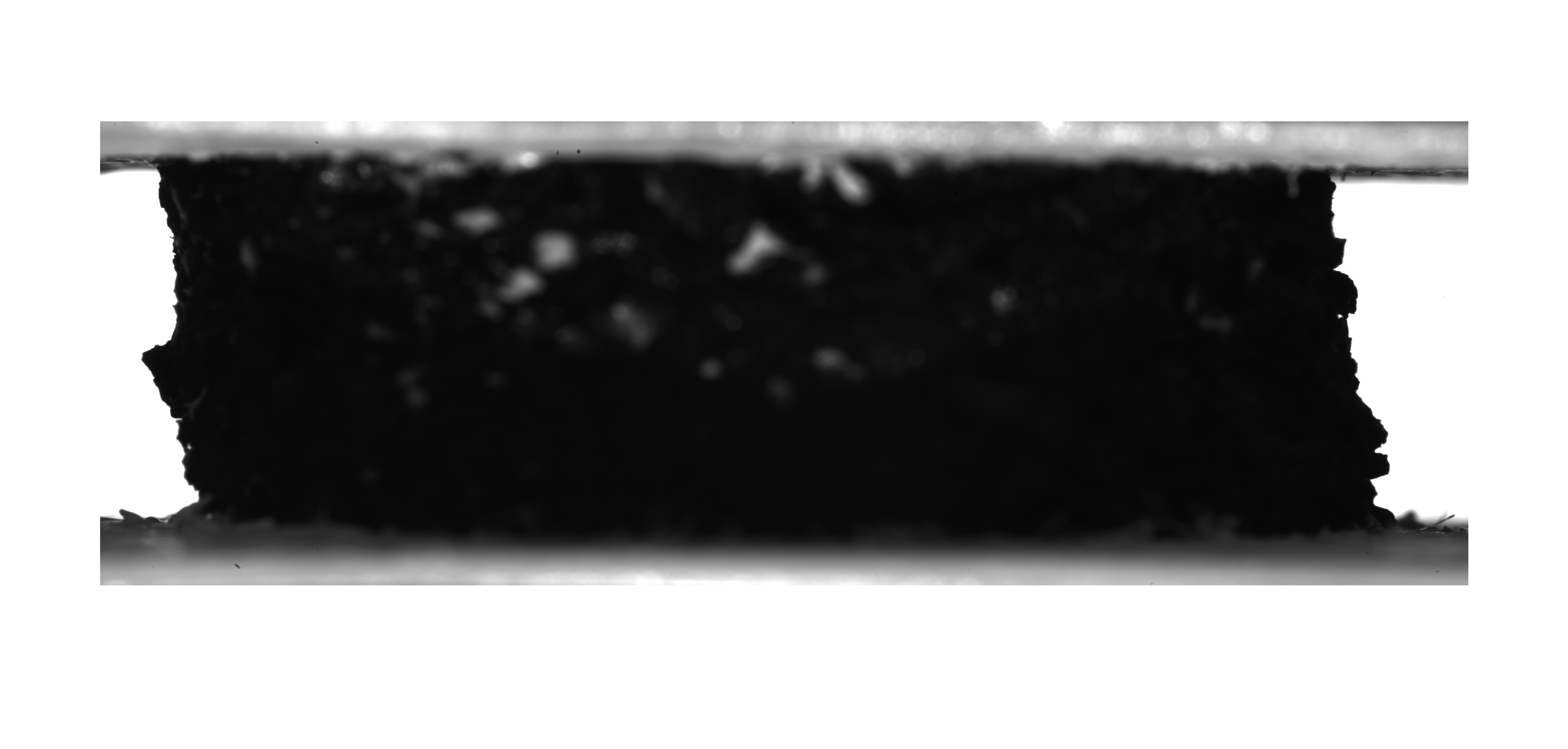}}
        \subfigure{\includegraphics[width=40mm,keepaspectratio]{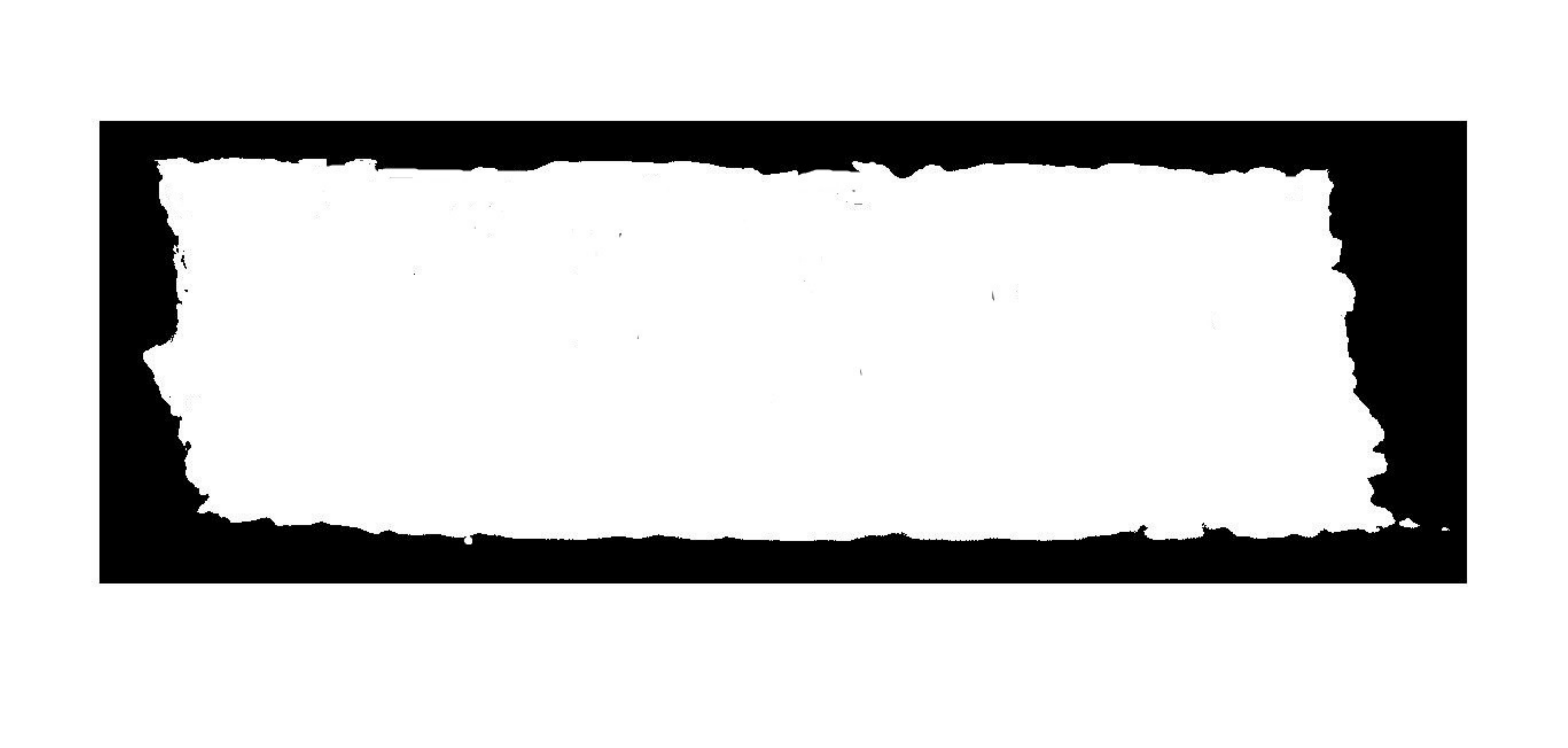}}\\
        \subfigure{70\%}{\includegraphics[width=40mm,keepaspectratio]{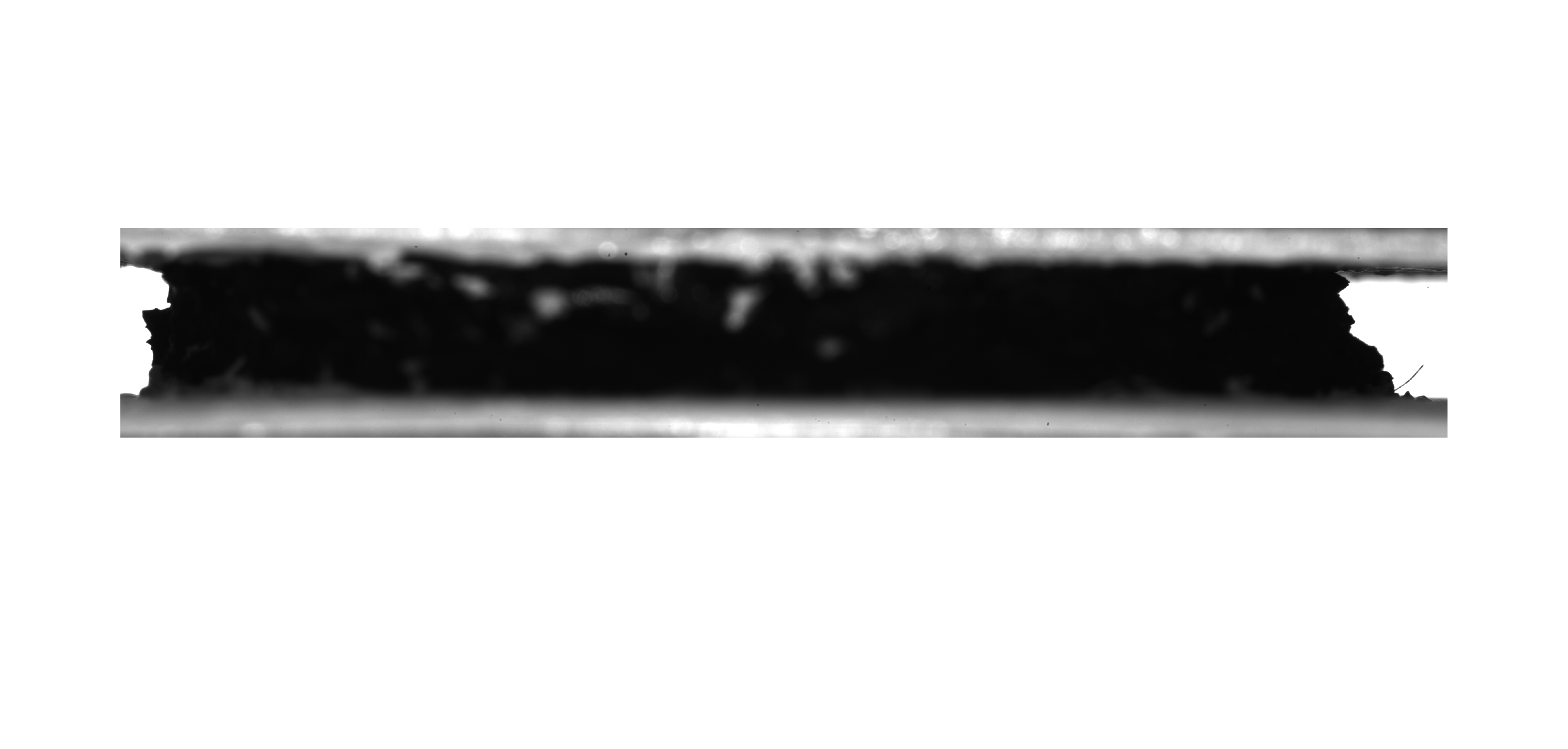}}
        \subfigure{\includegraphics[width=40mm,keepaspectratio]{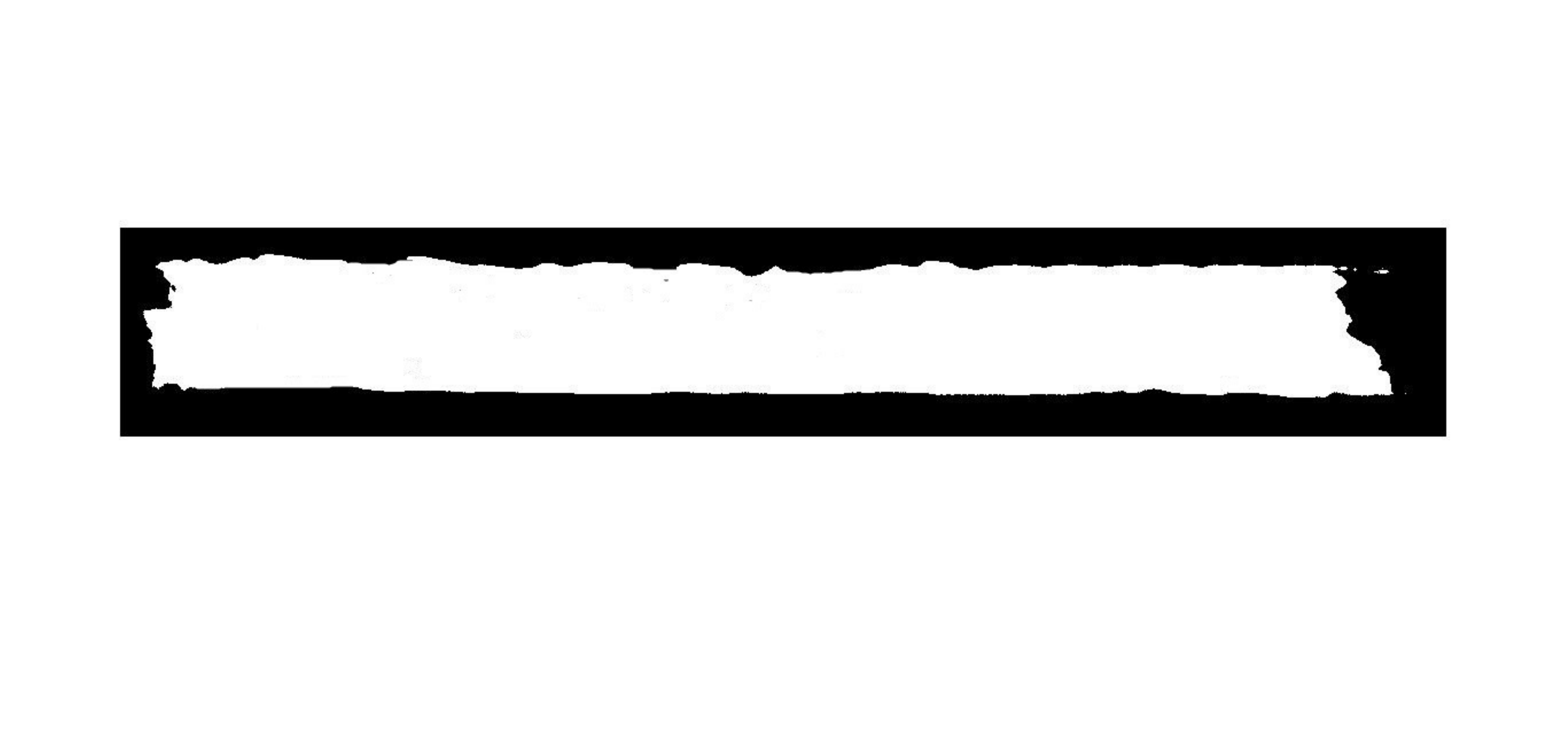}}
        \caption[Original image from the digital camera vs Image after processing ]{Original image from the digital camera (left column) versus image after processing (right column) for TDM 600 at 0\%, 35\%, 70\% defomation.}\label{fig:fig_pics_comp}
        \end{figure}

\begin{figure}[H] 
       \centering
       \includegraphics{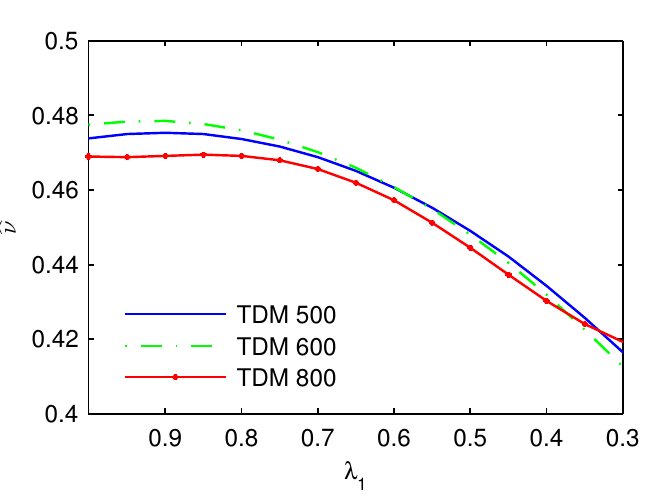} 
       \caption[Non-linear Poisson's coefficient]{Non-linear Poisson's coefficient $\widehat{\nu}$ evaluated during compression tests.}\label{fig:fig_nlpoisson}
       \end{figure}

Uniaxial compression tests were performed using a multi-step relaxation procedure. Thin Teflon sheets with lubricant were placed between platens and specimen surfaces. Specimen were cylinders with diameter of $\approx$ 27 mm and length of $\approx$ 15 mm die-cut from a sheet stock. At each step of the loading process, the strain level is increased by $5 \%$ at a strain rate of $\epsilon =0.01 \,\mathrm{s^{-1}}$ up to $70 \%$ strain. Between each loading step there is a $ 600$ s dwell to allow for relaxation of the material (Figure \ref{fig:fig_load}). We consider the value of the stress at the end of each dwell-interval as the equilibrium stress; these are shown as the red curves in Figs.\ \ref{fig:fig_load}(b)-\ref{fig:fig_load}(d).  Note that we take the $1$-direction to be the axis of compression.
In order to compare the model to the compression data, one requires knowledge of the material's three-dimensional state of deformation. Since the TDMs are compressible we can not make the usual rubber elasticity assumptions and require information on the materials transverse response. To evaluate the transverse behavior, pictures were taken at the end of every relaxation period with a digital camera mounted on a tripod; see Figure \ref{fig:fig_pics_comp}. The digital images were processed using the image processing toolbox in MATLAB \cite{Matlabi}.
The white area $ A $, measured in pixels in Figure \ref{fig:fig_pics_comp}, is the region occupied by the TDM sample, and was determined for each image. The height $ h $ of the sample is known at each step of the test  from the experimental controller. TDMs due to their composition are very difficult to cut and usually the samples do not have a straight edge. For this reason, we compute an average width of the sample as: 
    \begin{equation}
    \tilde{w}=\dfrac{A}{h}
    \end{equation}
and the average lateral strain as:
    \begin{equation}
    \widetilde{\log\lambda}_{2}=\widetilde{\log\lambda}_{3}=\log\dfrac{\tilde{w}_n}{\tilde{w}_{0}}\,,
    \end{equation} 
where $ \tilde{w}_n $ is the lateral dimension at step $ n $ and $ \tilde{w}_{0} $ is the lateral dimension at the undeformed state. In order to appreciate the non-linear compressibility of TDMs, we define, among several possibilities \cite{kakavas2000prediction,smith1999interpretation,Helf}, the non-linear Poisson's coefficient $\widehat{\nu}$ as the negative ratio of the lateral extension and axial contraction measured in the logarithmic strain:
    \begin{equation}
    \widehat{\nu}\colonequals-\dfrac{\log\lambda_2}{\log\lambda_1 }\, .
    \end{equation}
The measured values of $\hat\nu$ are shown in Fig.~\ref{fig:fig_nlpoisson} and display the material's distinct nonlinearity.\par
\begin{figure}[H]
       \subfigure[\,TDM 500]{\includegraphics{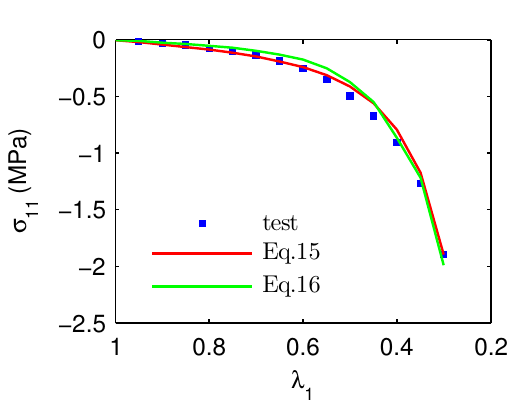}}
       \subfigure[\,TDM 600]{\includegraphics{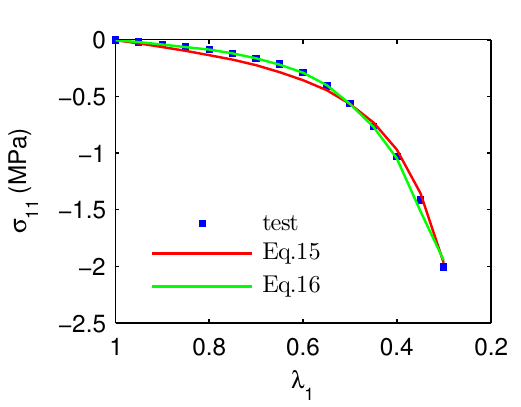}}
       \subfigure[\,TDM 800]{\includegraphics{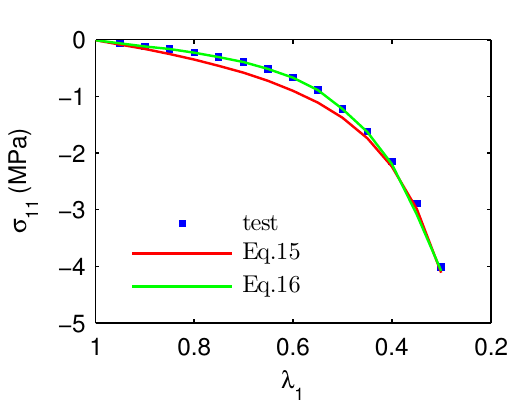}}
       \caption[Comparison between tests and exponentiated-hencky strain energy function]{Comparison between compression stress corresponding to modified exponentiated Hencky energy $W_{\mathrm{eHm}}$, equation \eqref{s} and \eqref{s2}, and experimental tests for different densities.} \label{fig:fig_compression} 
       \end{figure}

From equation \eqref{sigma}, if we consider \textit{s} the value of the uniaxial Cauchy stress, by projection on the Lie-algebra $\mathfrak{sl}$(\textit{n}) of trace-free tensors, we have
    \begin{equation}\label{taucomp}
    2\,\mu\, e^{k\,\|\dev_3\log {U}\|^{2}-\mathrm{tr}(\log {U})}\dev_3\log {U}= \dev_3\sigma=
    \left({\begin{array}{ccc}
    \dfrac{2}{3}s& 0 & 0 \\
    0 & -\dfrac{1}{3}s & 0\\
     0 & 0 & -\dfrac{1}{3}s
    \end{array} }\right)\,,
    \end{equation}
which leads to the requirement that under uniaxial stress, $U$ has the following form \cite{vallee1978lois}: 
    \begin{equation}\label{asatzU}
    U=\left({\begin{array}{ccc}
    e^{a+\frac{1}{3}x}& 0 & 0 \\
    0 & e^{-\frac{1}{2}a+\frac{1}{3}x} & 0\\
    0 & 0 & e^{-\frac{1}{2}a+\frac{1}{3}x}
    \end{array} }\right) = e^{\frac{1}{3}x}\left({\begin{array}{ccc}
    e^{a}& 0 & 0 \\
    0 & e^{-\frac{1}{2}a} & 0\\
    0 & 0 & e^{-\frac{1}{2}a}
    \end{array} }\right)\,.
    \end{equation}
This in turn leads to the stress expression
    \begin{equation}\label{s}
    s=3\,\mu\, e^{k\,\frac{3}{2}\,a^{2}-x}a\,,
    \end{equation}
where $a = \frac{2}{3}( \log\lambda_1-\widetilde{\log\lambda_2})$ and $x = \log\lambda_1+2\widetilde{\log\lambda_2}$ are the experimentally known measures of the deformation.
Using the measured values of $a$ and $x$ we can compare the model's predicted stress response from equation \eqref{s} to the measured stresses.  
Note that one can also project onto the spherical part of equation \eqref{sigma}, which gives 
    \begin{equation}\label{s2}
    s=3\left[
    \kappa\, e^{\hat{k}x^2}x+\kappa_1
    e^{\tilde{k}|x|^m}\frac{|x|^m}{x}\right]e^{-x}.
    \end{equation}
The comparison of equation \eqref{s} and \eqref{s2} to the experimental data is shown in Fig.~\ref{fig:fig_compression}. The results indicate good correlation. Note that the material parameter set is the same as used in the comparison of the shear data.
       
 \subsection{Pseudo-hydrostatic compression}
 
\begin{figure}[H]
       \subfigure[\,TDM 500]{\includegraphics{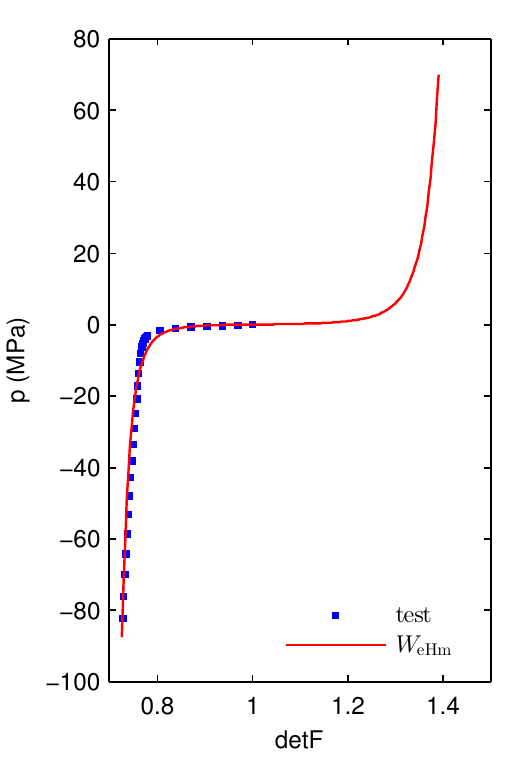}}
       \subfigure[\,TDM 600]{\includegraphics{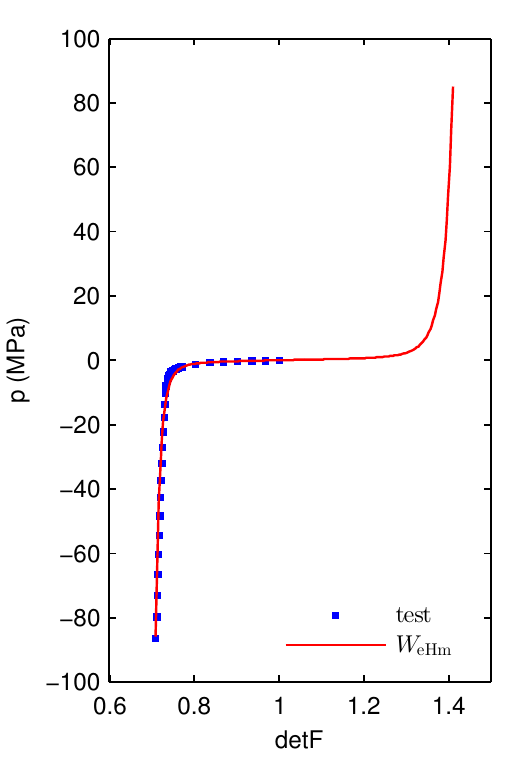}}
       \subfigure[\,TDM 800]{\includegraphics{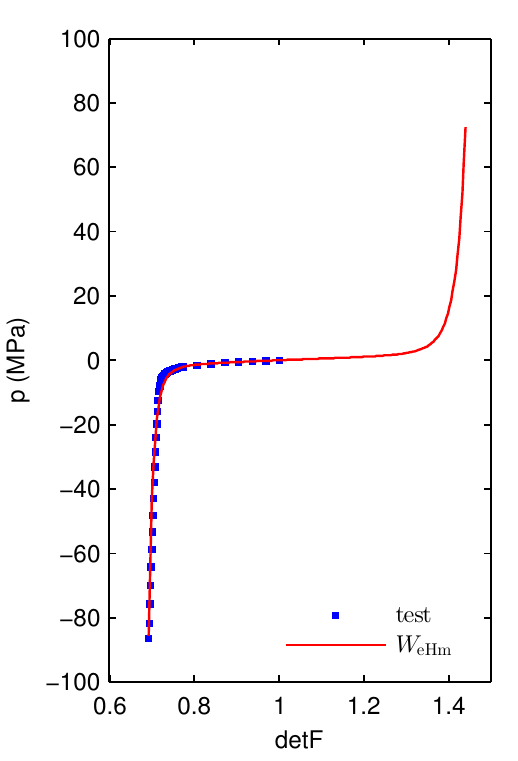}}
       \caption[Comparison between tests and exponentiated-hencky strain energy function]{Comparison between pseudo-volumetric response corresponding to equation \eqref{sigma} and experimental tests for different densities with $p \approx \sigma_{11}$.}\label{fig:fig_volumetric} 
       \end{figure}

As a third mode of deformation we consider an experiment that for quasi-incompressible  materials provides an approximation to a hydrostatic compression test.  Here we consider a lubricated cylindrical specimen that is inserted into a rigid (steel) cavity of the same radius and then axially compressed. During the test a force was applied on top of a steel piston at a volume ratio rate $0.0067\,\mathrm{s^{-1}}$. The procedure is explained in details in \cite{Montella2014}. Considering the axis of compression to be the $1$-direction, to good approximation this test
follows the kinematic path
    \begin{equation}
    \log U = \left( \begin{array}{ccc} \log\lambda_1 & 0 & 0\\ 0 & 0 & 0 \\ 0 & 0 & 0 \end{array}\right)\, ,
    \end{equation}
where in the experiment $\log\lambda_1$ is measured.  Likewise in the experiment $\sigma_{11} = s_{11} + p$ is measured,
where $s_{11}$ is the axial deviatoric stress component and $p$ is the pressure.  The experiment is designed to
test the pressure-volume relation.
The Jacobian of the
deformation, $\det F = \lambda_1$, is readily available from the experiment.  However the pressure is
approximated as $p \approx \sigma_{11}$, which is only valid for $p \gg s_{11}$.  For the present model, under
the given deformation state,
    \begin{equation}\label{p}
    \frac{s_{11}}{p} = \frac{2\mu\, e^{k\frac{2}{3}(\log\lambda_1)^2 -\log\lambda_1}\log\lambda_1}
    {\kappa\, e^{\hat{k}(\log\lambda_1)^2 -\log\lambda_1}\log\lambda_1
    +\kappa_1 e^{\tilde{k}|\log\lambda_1|^m -\log\lambda_1}\frac{|\log\lambda_1|^m}{\log\lambda_1} }\, .
    \end{equation}
For TDMs in this experiment, this ratio is not small enough at initial levels of deformation to result in a
valid pressure-volume experiment.  For example, for TDM 500, one must have $\det F \not\in (0.85,1.2)$ for the ratio
to take on values of less than $0.05$.  The plot of the data and the model prediction are shown in
Figure \ref{fig:fig_volumetric} and show good agreement.  Note that the plotted
pressure is approximated as $\sigma_{11}$ in both the model and the experiment for full consistency. It is to be noted that our exponentiated Hencky energy is crucial here. The quadratic Hencky energy leads to a pressure-volume relation that is not even invertible for $\det F > e$ \cite{vallee1978lois}. The exponentiated form itself alleviates this problem. The modified spherical term which we added to the original exponentiated Hencky energy allows for the sharp kink in
the pressure-volume relation due to void collapse.  As with the prior deformation modes, the results shown are produced with the \emph{exact same set of parameters}.
    
\section{Parameter identification equilibrium response}\label{sec:III}

For the three states of deformation under consideration, we have utilized a single set of parameters per material density. The parameter estimation itself was performed using the non-linear least square (NLS) optimization method to minimize the residuals. The main difficulty in such a NLS problem is to find a unique set of optimal parameters. Several numerical algorithms have been used in the literature to solve NLS problems \cite{bjorck1996,Aster}; they are usually a modification of the Newton method and require an initial guess for the solution. The iterative technique furnishes an optimal solution when some stopping criteria are met. In this paper we modified the function \textit{lsqcurvefit} in the optimization Toolbox of MATLAB \cite{Matlab} to fit the different sets of data together \cite{Ogden2004}. The initial guess used was the physical estimate of the moduli obtained by previous experiments on the TDMs (Table \ref{tab:II}); the remaining parameters were initialized to $k = \hat{k} = \tilde{k} = 0$ and $m=2$. We imposed positive values as lower bounds on the parameters; moreover we respect the mathematical bounds on the parameters \cite{NeffGhibaLankeit}. Table \ref{tab:III} gives the resulting optimized values which were used for modelling the behavior displayed in the prior sections. It is useful to note that the fit values for the shear, low strain bulk, and large strain bulk moduli are all sensibly close to the original values derived directly from the tangents to the experimental response.  The only oddity in values occurs for $\kappa_1$, the high compression bulk contribution which drops in going from TDM 600 to TDM 800, whereas an increase would intuitively be expected.  The cause for this is as yet unexplained but could be related to some sort of local material failure.  It should be noted that TDM's, due to their very porous structure, possess structural level failure mechanisms that initiate before the usually expected ones, such as filler dewetting \cite{bueche:60*1,bueche:61,rigbi:80} etc.\ as is common in filled elastomeric composites.
    \begin{table}[H]
    \caption{Initial guess for parameter identification procedure  from tangents to the experimental data$^*$.}\label{tab:II}
    \centering
    \begin{tabular}{C{2cm}C{1cm}C{1cm}C{1cm}}
    \toprule
    \multirow{2}{*}{Material} & $\mu$ & $\kappa$ & $\kappa_{1}$ \\
    & (MPa) &(MPa)& (MPa)\tablebreak
    \hline\\[-.7em]
    TDM $500$ & $0.22$ & $2.40$ &  $297$  \tablebreak
    TDM $600$ & $0.31$ & $2.70$ &  $315$  \tablebreak
    TDM $800$ & $0.63$ & $4.50$ & $281$  \tablebreak
    \hline\\[-.9em]
    \multicolumn{4}{l}{\textsuperscript{*}\footnotesize{Parameters obtained from \cite{Montella2014}.}}
    \end{tabular}
    \end{table}
    
    \begin{table}[H]
    \caption{Parameters - exponentiated-Hencky energy function.}\label{tab:III}
    \centering
    \begin{tabular}{C{2cm}C{1cm}C{1cm}C{1cm}C{1cm}C{1cm}C{1cm}C{1cm}}
    \toprule
    \multicolumn{8}{c}{\vphantom{\rule[-1.05em]{0em}{2.45em}}$\widehat{W}_{\mathrm{eHm}}(U)\colonequals
    \dfrac{\mu}{k}\, e^{k\,||\dev_3\log U||^{2}} +\dfrac{\kappa}{2\,\widehat{k}}\,e^{\widehat{k}\,[\tr(\log U)]^{2}}+\dfrac{\kappa_{1}}{m\,\tilde{k}}\,e^{\tilde{k}\,\left |\mathrm{tr}(\mathrm{\log U})\right|^{m}}$}\\[.5em]
    \hline\\[-.7em]
    \multirow{2}{*}{Material}{\vphantom{\rule{0em}{1.12em}}} & $\mu$ & $k$  &$\kappa$ & $\widehat{k}$ &$\kappa_{1}$ & $\tilde{k}$& $m$\\
    & (MPa) & (-) & (MPa)& (-) & (MPa)& (-)& (-)\tablebreak
    \hline\\[-.7em]
    TDM $500$ & $0.12$ & $0.59$ & $1.40$ & $0.13$& $116$ & $ 268 $ & $ 4 $\tablebreak
    TDM $600$ & $0.19$ & $0.39$ & $2.80$ & $0.13$  & $647$ & $ 1989 $ & $ 6 $\tablebreak
    TDM $800$ & $0.50$ & $0.27$ & $4.40$ & $0.13$  & $404$ & $ 1353 $ & $ 6 $\tablebreak
    \hline
    \end{tabular}
    \end{table}
%
%
\section{Non-equilibrium response}

\begin{figure}[H] 
        \centering
        \includegraphics{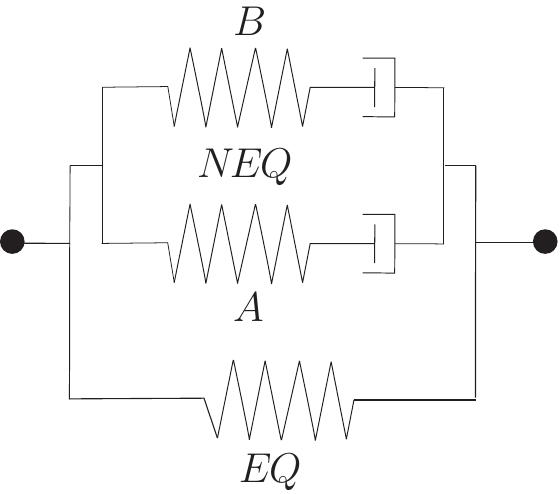} 
        \caption[Rheological model]{One dimensional rheological model for rate-dependent behavior of TDM.} \label{fig:fig_reol}
        \end{figure}

In experimental investigations, filler-reinforced rubber like SBR, the main component of TDMs, shows many nonlinear effects when subjected to dynamic loads. The main ones being the pronounced dependence of the material behaviour on the dynamic strain amplitude together with rate dependent response. The first, often termed the Payne-effect \cite{payne}, can be described as a reversible softening with increasing dynamic strain amplitude. 
To the authors' knowledge there are still no well-accepted models that incorporate both the Payne-effect and rate dependency.
Thus our goal in modeling the non-equilibrium behavior will be limited to the rate dependency of the material at fixed frequencies and amplitudes.
To that modest end, a finite strain model of viscoelasticity is constructed considering the multiplicative decomposition of the deformation gradient $F$ into elastic $F_{e}$ and inelastic $F_{i}$ parts as proposed by Sidoroff \cite{Sidoroff}. Here we assume the existence of two viscous mechanisms associated to the material: intermolecular resistance at the microscale level and grain interactions at the macroscale level. The first is associated with a Maxwell element including a non-linear spring (A) while the second is associated to a Maxwell element in which a linear spring is included (B). The choice of modeling the interaction between the rubber particle inside the TDM with a linear law is due to the presence of the binder at the grain interface. The binder acts as an internal constraint allowing only normal contact interaction between the grains. For this reason, both relative rotation and sliding, which are usually found in granular materials \cite{Oda,Jiang}, are not allowed or are negligible between the grains and fibers.
A one dimensional rheological schematic is presented in Figure \ref{fig:fig_reol}. 
If we consider for each viscous mechanism a set of internal variables $F_{i}^{k}\ (k=A,B)$ that can be viewed as the deformation gradient associated with each dashpot, then we can write the total free energy as:
    \begin{equation}\label{functional}
    W^\mathrm{v}_{\mathrm{eHm}}=W_\mathrm{eHm}^{\mathrm{EQ}} +W_{\mathrm{ NEQ}}^{A}(b_{e}^{A})+W_{\mathrm{NEQ}}^{B}(b_{e}^{B})\, ,
    \end{equation}
where $W_\mathrm{eHm}^{EQ}$ represents the strain energy in the equilibrium spring discussed in the previous part of this paper and $W_{NEQ}^{k}$ the strain energy in each Maxwell element associated to the ``elastic" left Cauchy deformation tensor $b_{e}^{k}= F_{e}^{k}\cdot [F_{e}^{k}]^{T}$, also called the Finger tensor.
For most polymer based materials, the volumetric deformation is purely elastic and the viscous effects are restricted to the isochoric component of the deformation. Following this assumption the strain energy for the Maxwell elements can be written as:
    \begin{align}\label{funcs}
    W_{\mathrm{NEQ}}^{A}(b_{e}^{A})= \dfrac{\mu_A}{k_A}\, e^{k_A\,||\dev_3\log b_{e}^{A}||^{2}}\\
    W_{\mathrm{NEQ}}^{B}(b_{e}^{B})= \mu_B\,\,||\dev_3\log b_{e}^{B}||^{2}\, .
    \end{align}
The general theory of viscoelasticity at finite strains used in this work follows the developments of \cite{Govindjee1998,govi1997}. Here we recall only the essential equations.
As a consequence of the Clausius-Duhem inequality, the Kirchhoff stress is given as
    \begin{equation}
    \tau = \tau_{EQ} + \sum \tau_{NEQ}^k\, ,
    \end{equation}
where $\tau_{NEQ}^k = 2[\frac{\partial W^k_{NEQ}}{\partial b_e^k}] b_e^k$. Consistent with the Clausius-Duhem inequality, the evolution of $b_e^k$ is given by:
    \begin{equation}\label{evol_eq}
    \dfrac{1}{2}\,\mathcal{L}_{v}\,b_{e}^{k}\cdot[b_{e}^{k}]^{-1}=[\mathcal{V}^{k}]^{-1}:\tau_{\mathrm{NEQ}}^{k}\, ,
    \end{equation}
where $\mathcal{L}_{v}\,b_{e}^{k} = F \frac{\mathrm{d}}{\mathrm{dt}}[C_i^k]^{-1}F^T$ is the Lie derivative  of $b_e^k$ along the velocity field of the material motion, $C_i^k = [F_i^k]^T F_i^k$,  and $[\mathcal{V}^{k}]^{-1}$ is an isotropic fourth order fluidity tensor defined as:
    \begin{equation}\label{iso_tens}
    [\mathcal{V}^{k}]^{-1}=\dfrac{1}{2\,\eta_{D}^{k}}\left(\mathbb{1}^{4}-\dfrac{1}{3}\mathbb{1}\otimes\mathbb{1}\right).
    \end{equation}
Here $\mathbb{1}^{4}$ is the fourth order symmetric identity tensor, while $\eta_{D}^{k}>0$ represents the deviatoric viscosities. In our model $\eta_{D}^{A}=12 \,\mathrm{s\cdot N \cdot mm^{-2}}$ and $\eta_{D}^{B}=1 \,\mathrm{s\cdot N \cdot mm^{-2}}$ for all the different densities and for all the testing modes.
The model presented implicitly defines the total Kirchhoff stress $\tau$. The actual use of the model requires the solution of the nonlinear relation \eqref{evol_eq} which we perform using the predictor-corrector method advocated in \cite{govi1997} and \cite{Govindjee1998}.

\subsection{Dynamic shear test}

Dynamic shear tests were performed at TARRC using the dual lap set up with samples of the same dimension as used in the static tests. For each sample the displacement was driven up to 33\% and 100\% of the initial thickness and the tests, for each amplitude, were carried out for 10 cycles. Further, each test was conducted at two frequencies 0.1 Hz and 1 Hz.
    \begin{figure}[H]
        \subfigure[\,TDM 500]{\includegraphics{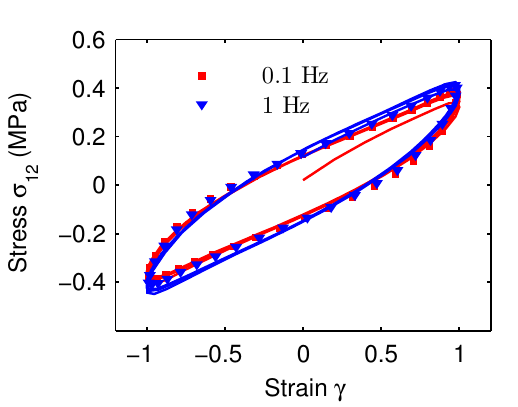}}
        \subfigure[\,TDM 600]{\includegraphics{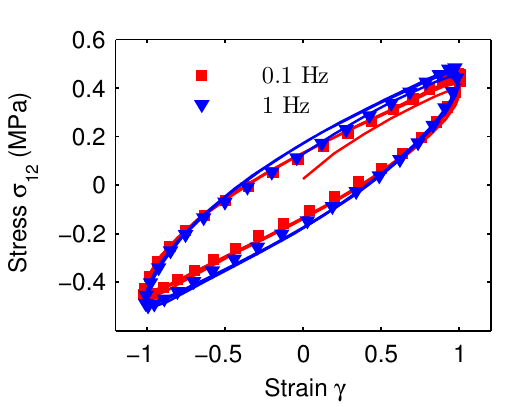}}
        \subfigure[\,TDM 800]{\includegraphics{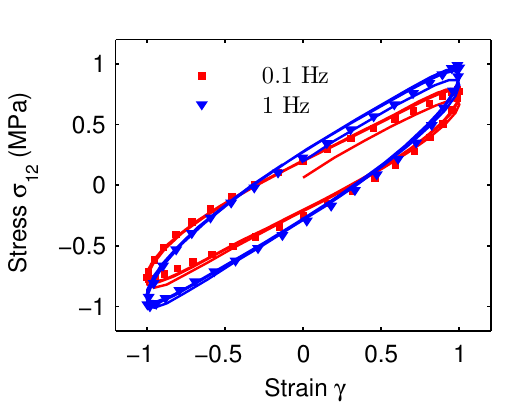}}
        \caption[Comparison between tests and exponentiated-Hencky strain energy function]{Comparison between cyclic shear tests (markers) and the viscoelastic model based on the modified exponentiated Hencky energy $W^\mathrm{v}_{\mathrm{eHm}}$ (solid line), equation \eqref{functional}, for different frequencies at 100\% amplitude.}\label{fig:fig_shear_dyin1}
        \end{figure}
        
        \begin{figure}[H]
        \subfigure[\,TDM 500]{\includegraphics{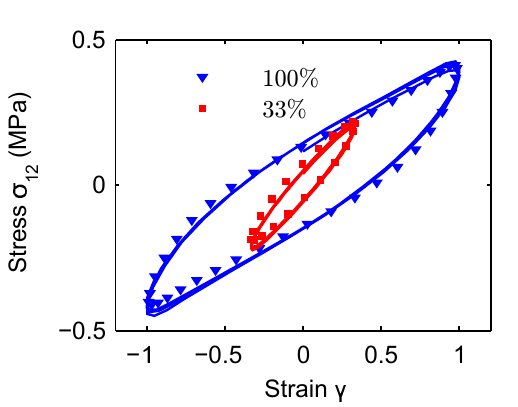}}
        \subfigure[\,TDM 600]{\includegraphics{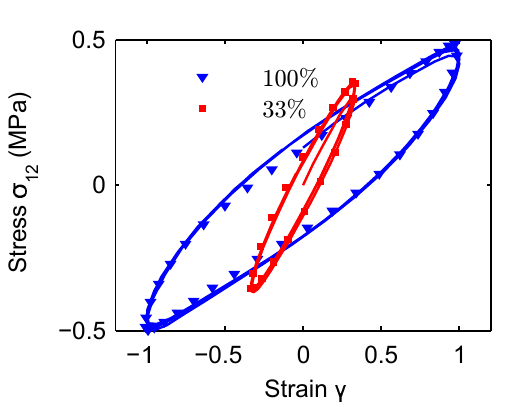}}
        \subfigure[\,TDM 800]{\includegraphics{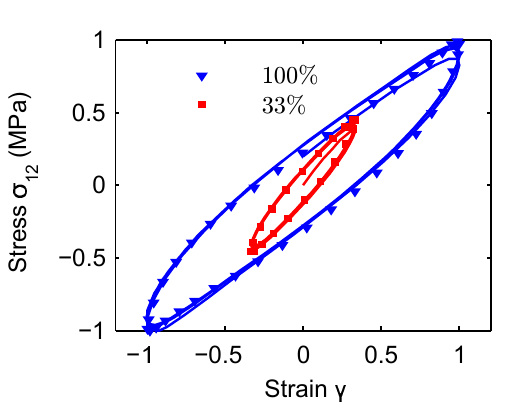}}
        \caption[Comparison between tests and exponentiated-Hencky strain energy function]{Comparison between cyclic shear tests (markers) and the viscoelastic model based on the modified exponentiated Hencky energy $W^\mathrm{v}_{\mathrm{eHm}}$ (solid line), equation \eqref{functional}, for different amplitudes at 1 Hz.}\label{fig:fig_shear_dyin2}
        \end{figure} 
Figure \ref{fig:fig_shear_dyin1} shows the dynamic response at 100\% strain amplitude for two frequencies.  Over this range of frequencies the material only weakly depends on the strain rate.  Figure \ref{fig:fig_shear_dyin2} considers 100\% and 33\% strain amplitude at a loading frequency of 1 Hz. Here one observes a strong amplitude dependent response.
Also shown in Figs.\ \ref{fig:fig_shear_dyin1} and \ref{fig:fig_shear_dyin2} are the predictions from fitting the model to the data.  The match is seen to be quite acceptable but it should be emphasized that due to the Payne-effect the values of the non-equilibrium parameters are frequency and amplitude dependent, as discussed more fully below.

\subsection{Dynamic compression test}

\begin{figure}[H]
        \subfigure[\,TDM 500]{\includegraphics{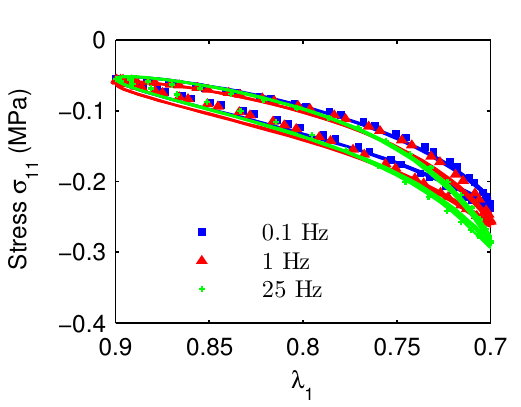}}
        \subfigure[\,TDM 600]{\includegraphics{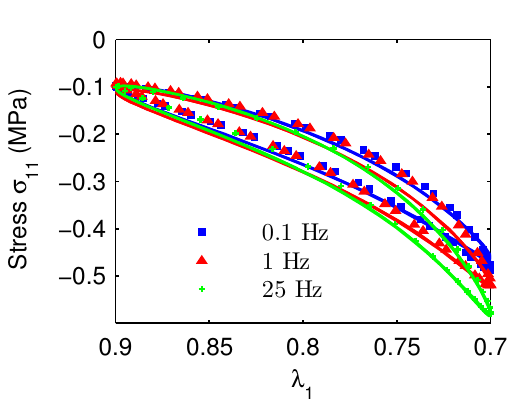}}
        \subfigure[\,TDM 800]{\includegraphics{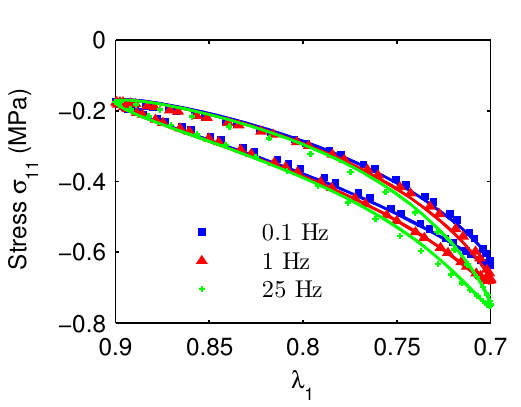}}
        \caption[Comparison between tests and exponentiated-Hencky strain energy function]{Comparison between cyclic compression tests (markers) and the viscoelastic model based on the modified exponentiated Hencky energy $W^\mathrm{v}_{\mathrm{eHm}}$ (solid line), equation \eqref{functional}, for different frequencies at 20\% amplitude.}\label{fig:fig_comp_dyin1}
        \end{figure}

\begin{figure}[H]
        \subfigure[\,TDM 500]{\includegraphics{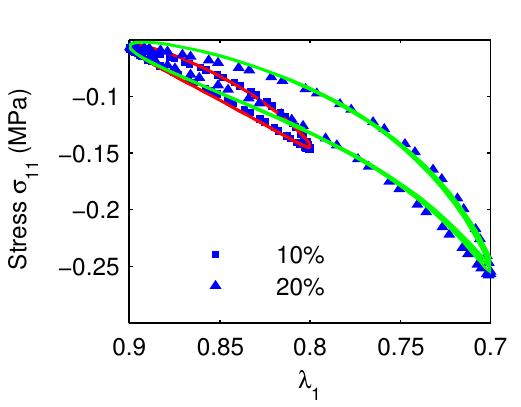}}
        \subfigure[\,TDM 600]{\includegraphics{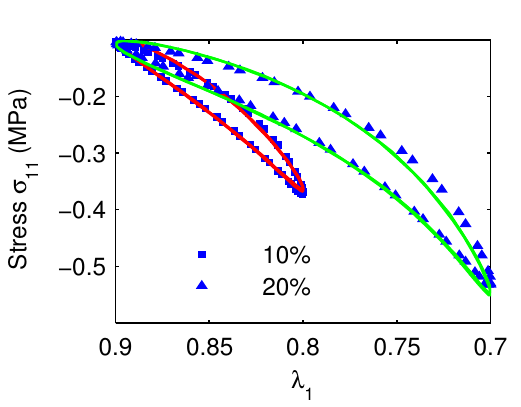}}
        \subfigure[\,TDM 800]{\includegraphics{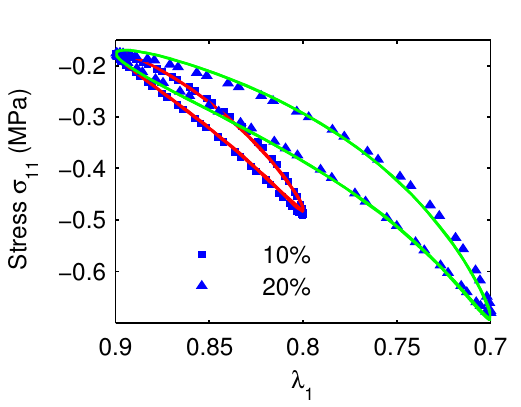}}
        \caption[Comparison between tests and exponentiated-Hencky strain energy function]{Comparison between cyclic compression tests (markers) and the viscoelastic model based on the modified exponentiated Hencky energy $W^\mathrm{v}_{\mathrm{eHm}}$ (solid line), equation \eqref{functional}, for different amplitudes at 1 Hz.}\label{fig:fig_comp_dyin2}
        \end{figure}

Uniaxial compression tests were carried out on a Bose Electroforce machine in a frequency range of 0.1 Hz to 25 Hz with the same setup as the static tests. The strain history consists of a static pre-strain of 10\% and a superimposed sinusoidal excitation varying in amplitude in the range of 1\% to 20\%.
Figures \ref{fig:fig_comp_dyin1} and \ref{fig:fig_comp_dyin2} show a few tests (data are shown as markers) from the many performed as they are representative of the overall abilities of the model. The model parameters used to generate the solid lines are discussed in the next section.
Figure \ref{fig:fig_comp_dyin1} shows the steady-state hysteresis curves at frequencies 0.1 Hz, 1 Hz and 25 Hz
at a constant strain amplitude of 20\%. It indicates that the stress increases with increasing frequencies and the material is stiffer at higher frequency. The correlation between model and experiment is seen to be good.
Figure \ref{fig:fig_comp_dyin2} shows the steady-state hysteresis curves with two different strain amplitudes, 10\% and 20\%, at 1 Hz frequency. These graphs confirm that the material subjected to smaller compressive  strain amplitudes is stiffer than material subjected to larger compressive strain amplitudes, similar to what was seen in the shear tests and the correlation between data and experiment is good. Other specimens showed a similar pattern even though they were taken about different mean strains and frequencies. They are not shown in this paper for brevity.
As part of the compressive strain campaign, we also evaluated the energy dissipated per hysteresis cycle as:
    \begin{equation}
    D=\oint \!\sigma_{11} d\lambda_1\, .
    \end{equation}

     \begin{figure}[H]
        \subfigure[\,TDM 500]{\includegraphics{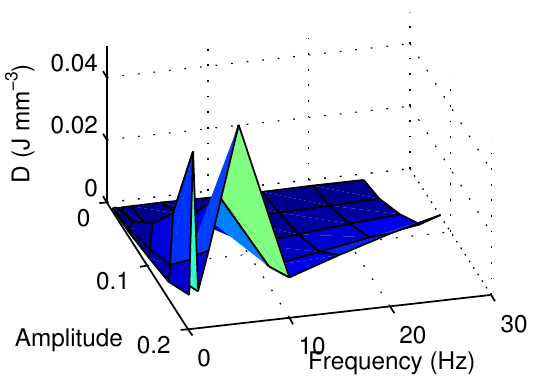}}
        \subfigure[\,TDM 600]{\includegraphics{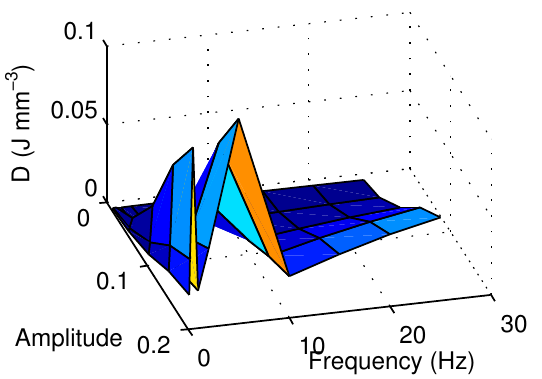}}
        \subfigure[\,TDM 800]{\includegraphics{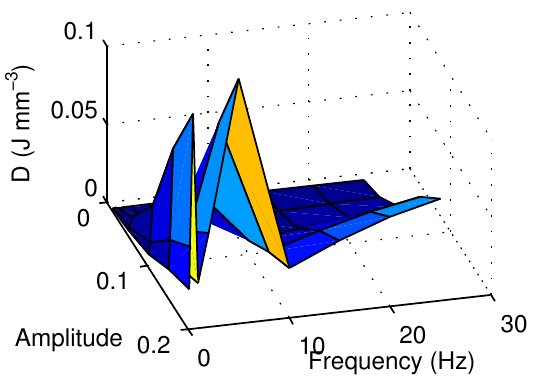}}
        \caption[Energy]{Energy dissipation per hysteresis cycle in compression.}\label{fig:fig_comp_diss}
        \end{figure}
\noindent The results, which are shown in Figure \ref{fig:fig_comp_diss}, display the presence of two transition regions which validates our use of two Maxwell elements in our model for TDMs in this range of amplitude and frequency. This in large part contributes to the good ability of the model to capture the hysterisis curves shown in Figs.\ \ref{fig:fig_shear_dyin1} - \ref{fig:fig_comp_dyin2}.

\section{Parameter identification viscous model}

\begin{figure}[H]
        \subfigure[\,TDM 500]{\includegraphics{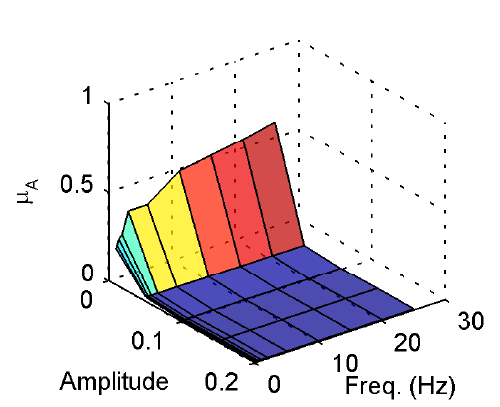}}
        \subfigure[\,TDM 600]{\includegraphics{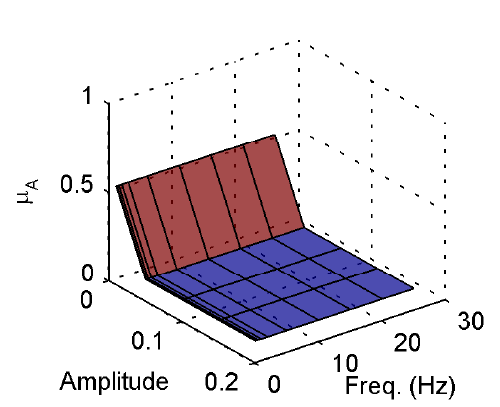}}
        \subfigure[\,TDM 800]{\includegraphics{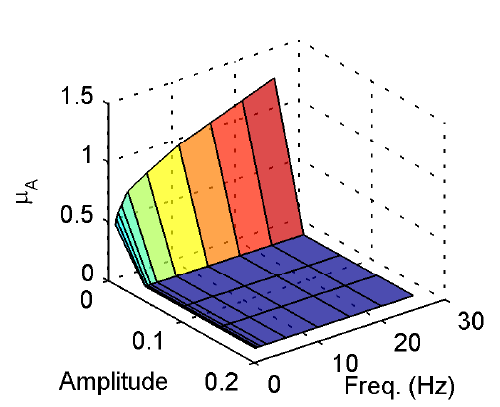}}
        \caption[Model Parameter]{Model parameter $\mu_A$ in the frequency and amplitude range.}\label{fig:fig_par_1a}
        \end{figure}

\begin{figure}[H]
        \subfigure[\,TDM 500]{\includegraphics{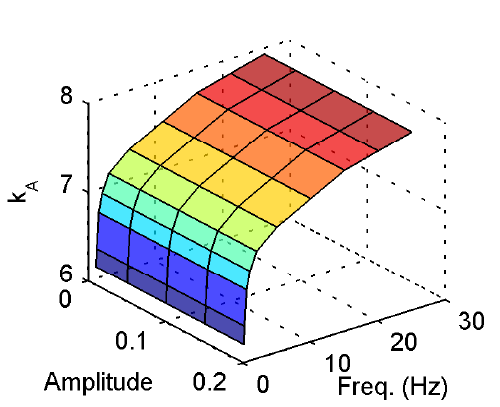}}
        \subfigure[\,TDM 600]{\includegraphics{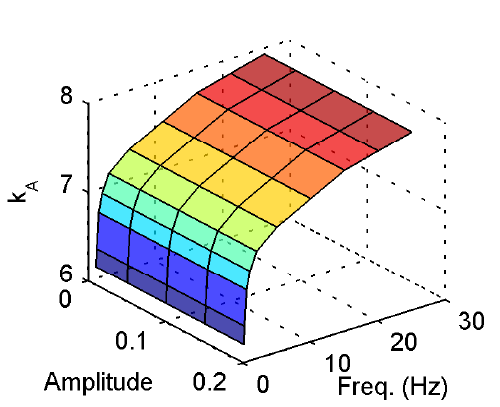}}
        \subfigure[\,TDM 800]{\includegraphics{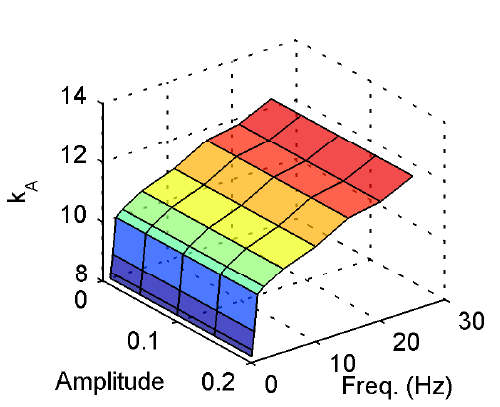}}
        \caption[Model Parameter]{Model parameter $k_A$ in the frequency and amplitude range.}\label{fig:fig_par_2a}
        \end{figure}

\begin{figure}[H]
        \subfigure[\,TDM 500]{\includegraphics{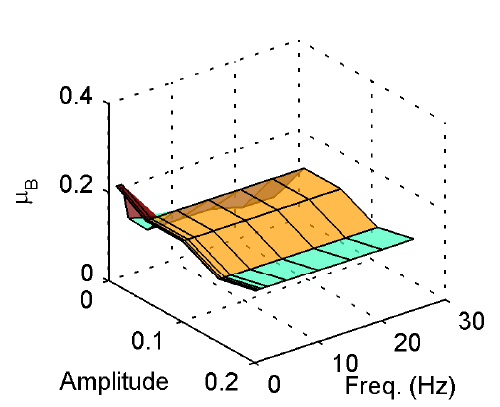}}
        \subfigure[\,TDM 600]{\includegraphics{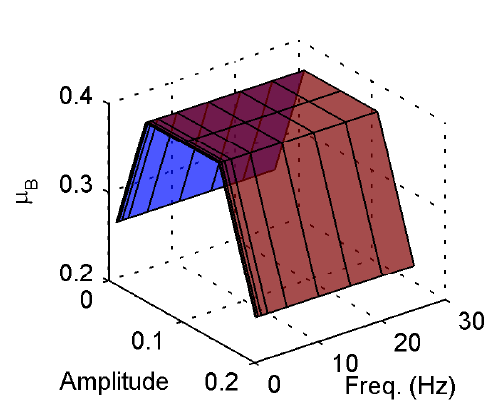}}
        \subfigure[\,TDM 800]{\includegraphics{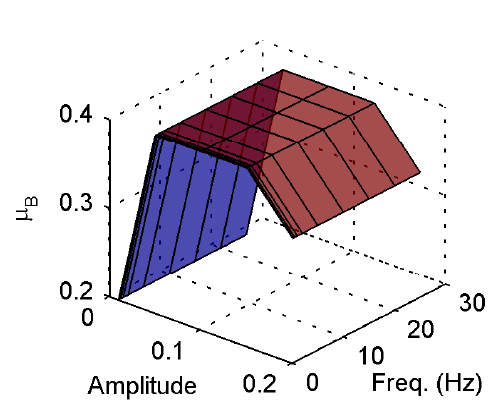}}
        \caption[Model Parameter]{Model parameter $\mu_B$ in the frequency and amplitude range.}\label{fig:fig_par_1b}
        \end{figure}

Since the response for large deformations is not a perfect sinusoid, the hysteresis cycles are not elliptical.
Therefore, the classical definition of storage and loss modulus is inapplicable. In this section we take a look at the parameters for each dissipation mechanism and show how they vary with frequency and amplitude.
The parameters associated with mechanism A, which we term the microscale level, show both amplitude and frequency dependence (Figure \ref{fig:fig_par_1a},\ref{fig:fig_par_2a}).
The parameter $\mu_A$ decreases with amplitude then it stays constant both in frequency and amplitude. It is mainly amplitude dependent reproducing the Payne-effect well known to be present in filler-reinforced rubber.  The parameter $k_A$ captures the frequency dependency of the material and it is constant with amplitude.
The single parameter associated with the mechanism B varies with the amplitude and stays constant with frequency (Figure \ref{fig:fig_par_1b}).

\section{Conclusion}

In this study we presented a hyper-visco-elastic constitutive model for TDMs to characterize the response of this class of materials under different deformation modes. The new model is based on an exponentiated Hencky strain energy that as shown in previous work improves several mathematical properties with respect to the classical quadratic Hencky energy function. There are two main advantages of the newly proposed model. The first one is its ability to describe different deformation modes with a unique set of parameters in the equilibrium range. The second advantage is that most of the parameters have a physical meaning simplifying the parameter fitting procedure.
An extensive experimental campaign on TDMs was conducted with both static and dynamic tests.
It was not the objective of this work to provide a model capable of describing dynamic characteristics of TDMs that includes simultaneous frequency and amplitude dependent effects, but rather to explore this new material and find the simplest model to characterize it for specific frequencies and amplitude.
The predicted results are in excellent agreement with the presented data and thus give a viable model for engineering applications of TDMs.

\section{Acknowledgments}
The work of G. Montella was carried out at the University of California, Berkeley under the grant ``Dottorato di Ricerca in Azienda” - POR Campania FSE 2007-2013, Asse IV.
The authors gratefully acknowledge Isolgomma s.r.l.\ for supplying the Tire Derived Material used in the tests and the Tun Abdul Razak Research Centre (TARRC) for providing the data for the shear tests, as well as I.-D.~Ghiba and R.\,J.~Martin (University of Duisburg-Essen) for their helpful remarks.

%
\bibliographystyle{plain} 
\bibliography{mybibfile}

\appendix
\section{Classical model fits}\label{sec:appA}
\begin{figure}[ht!]
        \subfigure[\,TDM 500]{\includegraphics{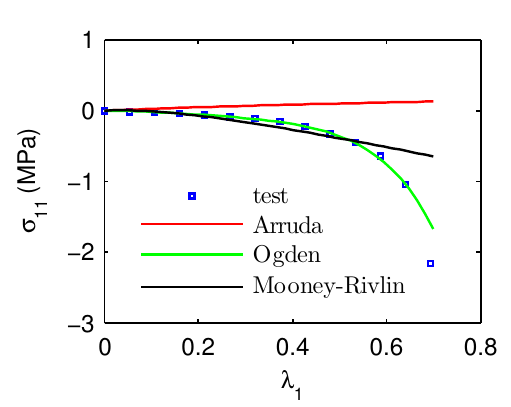}}
        \subfigure[\,TDM 600]{\includegraphics{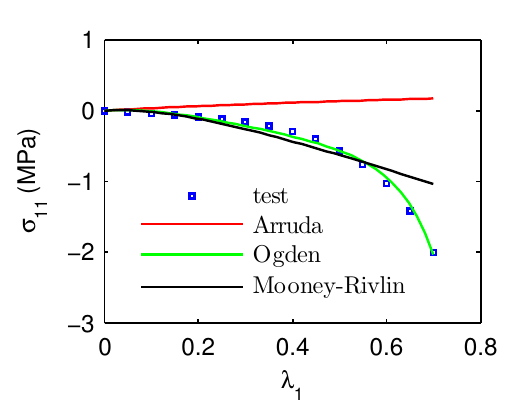}}
        \subfigure[\,TDM 800]{\includegraphics{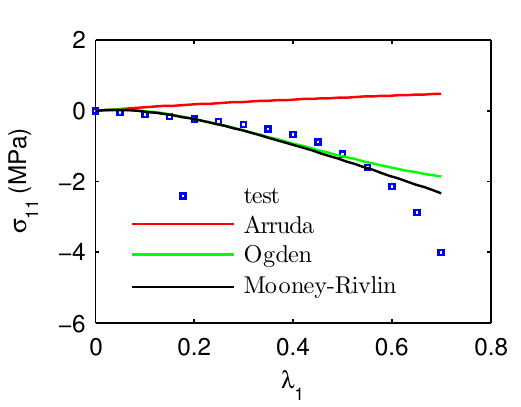}}
        \caption[Comparison between tests and strain energy functions]{Comparison between compression stress corresponding to different hyperelastic models and experimental tests for different densities.}\label{fig:fig_comp_other}
        \end{figure}

\begin{figure}[ht!]
        \subfigure[\,TDM 500]{\includegraphics{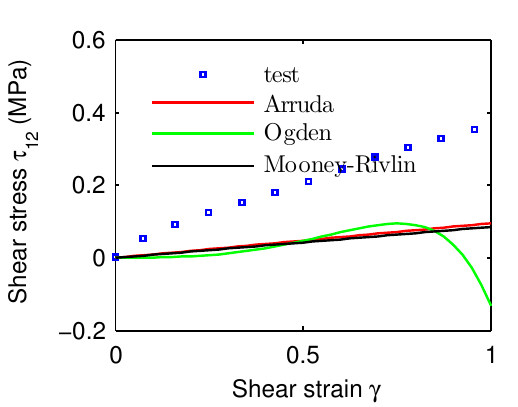}}
        \subfigure[\,TDM 600]{\includegraphics{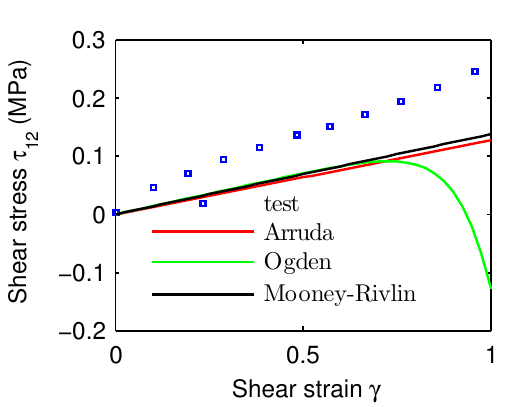}}
        \subfigure[\,TDM 800]{\includegraphics{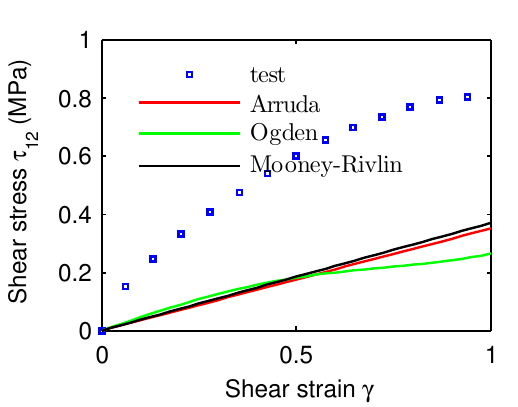}}
        \caption[Comparison between tests and strain energy functions]{Comparison between shear stress corresponding to different hyperelastic models and experimental tests for different densities.}\label{fig:fig_shear_other}
        \end{figure}

\begin{figure}[ht!]
        \subfigure[\,TDM 500]{\includegraphics{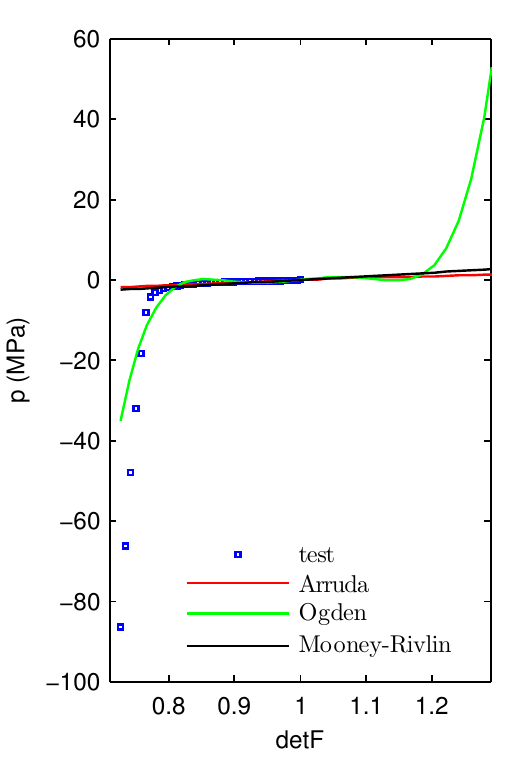}}
        \subfigure[\,TDM 600]{\includegraphics{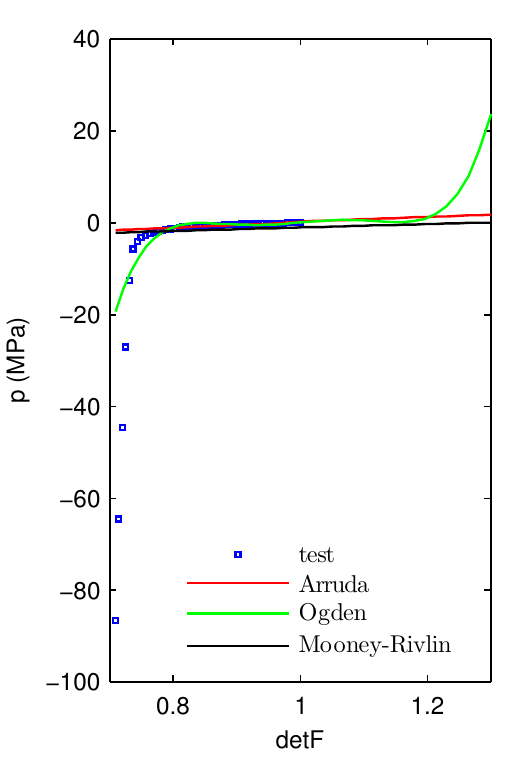}}
        \subfigure[\,TDM 800]{\includegraphics{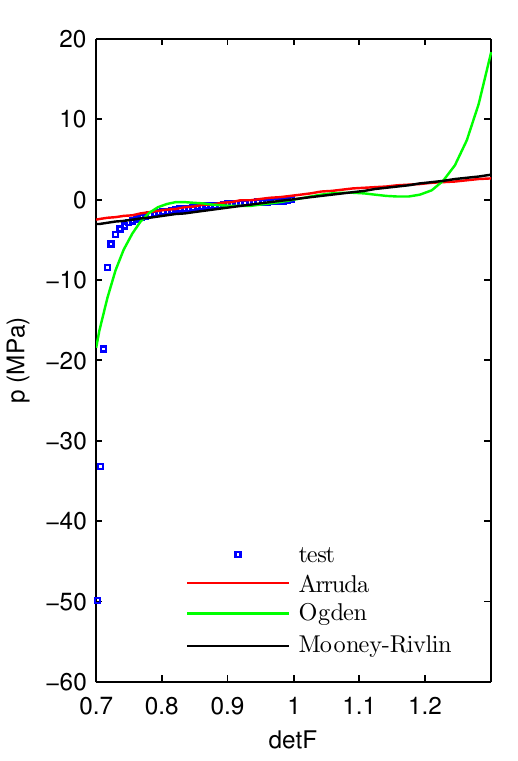}}
        \caption[Comparison between tests and strain energy functions]{Comparison between pseudo-volumetric response for different hyperelastic models and experimental tests for different densities.}\label{fig:fig_vol_other}
\end{figure}

Montella, Calabrese, and Serino \cite{Montella2014} attempted to fit classical hyperelastic models to TDM response but showed that the fits in general were poor.  With our additional data, we find similar (if not worse) results.  
Using the same fitting method as was done with our new exponentiated Hencky model, one has the following results when trying to utilize the Arruda-Boyce model \cite{Arruda1993}, the Mooney-Rivlin model \cite{mooney1940}, and the Ogden model truncated to the third order expansion \cite{ogden1972} for compression Fig.\ \ref{fig:fig_comp_other}, shear Fig.\ \ref{fig:fig_shear_other}, and pseudo-compression Fig.\ \ref{fig:fig_vol_other}.


\null
\newpage
\null
\newpage

\end{document}